\documentclass[10pt]{amsart}
\usepackage{a4wide}
\usepackage{amsmath,amssymb,amsfonts,amsthm}
\usepackage[notref,notcite]{showkeys}
\usepackage{mathtools}
\usepackage{xcolor}
\usepackage{xargs}
\usepackage{dutchcal}
\usepackage{gnuplot-lua-tikz}
\usepackage{tikz}
\usepackage[colorlinks,citecolor=blue,urlcolor=blue,bookmarks=false,hypertexnames=true]{hyperref} 

\theoremstyle{plain}

\newtheorem{theorem}{Theorem}[section]
\newtheorem{lemma}[theorem]{Lemma}
\newtheorem{proposition}[theorem]{Proposition}
\newtheorem{corollary}[theorem]{Corollary}
\theoremstyle{definition}
\begingroup
\newtheorem{remark}[theorem]{Remark}
\newtheorem{definition}[theorem]{Definition}
\endgroup

\numberwithin{equation}{section}

\def\Div{\textup{div}}

\def\R{\mathbb{R}}
\def\Q{\mathbb{Q}}
\def\Z{\mathbb{Z}}
\def\N{\mathbb{N}}
\def\dist{\textup{dist}}

\def\P{\mathcal P}

\newcommand{\po}{{\phi^\circ}}

\newcommand{\virg}[1]{``#1''}
\def\h{^\e}

\let\div\undefined 
\def\div{\textnormal{div}}
\newcommand{\ud}{\,\textnormal{d}}
\newcommand{\sd}{\,\textnormal{sd}}
\newcommand{\ez}{\e\Z^N}
\newcommand{\kz}{\e_k\Z^N}
\newcommand{\W}{\mathcal{W}}
\def\z{\mathbf{z}}

\newcommand{\e}{\varepsilon}

\newcommand{\pa}{\partial}
\newcommand{\beq}{\begin{equation}}
\newcommand{\eeq}{\end{equation}}
\newcommand{\medint}{-\kern -,375cm\int}

\usepackage[backend=bibtex, maxbibnames=99]{biblatex}

\title{Discrete-to-{continuum} crystalline curvature flows}

\author[A. Chambolle]{Antonin Chambolle}
\address[Antonin Chambolle]{CEREMADE, CNRS, Universit\'e Paris-Dauphine, PSL, and Mokaplan, INRIA Paris, France.}
\email[A. Chambolle]{antonin.chambolle@ceremade.dauphine.fr}

\author[D. De Gennaro]{Daniele De Gennaro}
\address[Daniele De Gennaro]{CEREMADE, CNRS, Universit\'e Paris-Dauphine, PSL, France, and Dipartimento di Scienze Matematiche Fisiche e Informatiche, Universit\`a di Parma, Italy.}
\email[D. De Gennaro]{degennaro@ceremade.dauphine.fr}
\author[M. Morini]{Massimiliano Morini}
\address[Massimiliano Morini]{Dipartimento di Scienze  Matematiche Fisiche e Informatiche, Universit\`a di Parma, Italy}
\email[M. Morini]{massimiliano.morini@unipr.it}

\keywords{Crystalline flows, dicrete-to-continuous approximation,
	 minimizing movements, distributional solutions}
\thanks{{\it Aknowledgements.}\! 
Antonin Chambolle is partially funded by a ``France 2030'' support managed by the Agence Nationale de la Recherche, under the reference ANR-23-PEIA-0004. Daniele De Gennaro has received funding from the European Union's Horizon 2020 research and innovation programme under the Marie Sk{\l}odowska-Curie grant agreement No~94532. Massimiliano Morini is a member of the Italian National Rsearch Group GNAMPA (INDAM) and is partially supported by PRIN 2022 Project ``Geometric Evolution Problems and Shape Optimization (GEPSO)'', PNRR Italia Domani, financed by European Union via the Program NextGenerationEU, CUP\_D53D23005820006. He wishes to warmly thank the hospitality of CEREMADE, where part of this research was conducted. The authors want to thank the anonymous  referees for the careful reading of the manuscript and their
comments, which helped improve the paper.
}
\begin{document}
	\maketitle
	\begin{abstract}
		We consider here a fully discrete variant of the implicit variational
		scheme for mean curvature flow~\cite{AlmTayWan,LucStu}, in a setting where
		the flow is governed by a crystalline surface tension defined by the limit of pairwise
		interactions energy on the discrete grid.
		The algorithm is based on a new discrete distance from the evolving sets, which prevents the occurrence of the spatial drift and pinning phenomena identified in \cite{MisiatsYip16,BraGelNov} in a similar discrete framework. We provide the first rigorous  convergence result  holding in any dimension, for any initial set and for a large class of purely crystalline anisotropies, in which the spatial discretization mesh can be of the same order or coarser than the time step.
	\end{abstract}
	
	\section{Introduction}\label{sec:intro}
	In this paper we analyse a space- and time-discrete approximation of crystalline mean curvature flows of the form
	\begin{equation}\label{evol law}
		V(x,t)=-\phi(\nu_{E(t)}(x))\kappa^\phi_{E(t)}(x), \qquad x\in\partial E(t), \quad t\ge 0,
	\end{equation}
	for a class of crystalline norms $\phi$. We recall that an anisotropy $\phi$ is said to be crystalline if and only if  $\{\phi\le 1\}$ is a polytope (or 
	equivalently, $\phi$ is the support function of a polytope). Moreover in the current paper we restrict ourselves to the case where $\{\phi\le 1\}$ is a 
	zonotope with rational generators~\cite{McM,ChaBra23}. 
	Here $V(x,t)$ stands for the (outer) normal velocity of the boundary $\partial E(t)$ at $x$, $\phi$ is a crystalline norm on $\R^N$ representing the surface tension, $\kappa^\phi_{E(t)}$ is the crystalline mean curvature of $\partial E(t)$ associated to  $\phi$, and $\nu_{E(t)}$ is the outer unit normal to $\partial E(t)$. The evolution law \eqref{evol law} has been considered to describe some phenomena in materials science and
	crystal growth; see e.g. \cite{Gurtin93,Taylor78}. Our main result is a convergence result of
	the discrete approximation to the continuous evolution, as the time and space steps go to zero,
	even in the somewhat surprising case where the space-step is greater or equal to the time-step.
	
	From the mathematical point of view, the lack of regularity of the differential operator involved in the definition of the crystalline curvature (see~\cite{BeNoPa,BePa96}) is the main reason why the  well-posedness of the crystalline mean curvature flow in every dimension has been a long-standing open problem.
	After some partial results (see for instance \cite{AlmTay95,AngGu89,BelCaChaNo,CasCha,GigaGiga,GigGigPoz14,GiGuMa98}), important breakthroughs have been obtained simultaneously in \cite{GigPoz16,GigPoz18,GigPoz20}, where a suitable crystalline theory of viscosity solutions was developed, and 
	with a different approach in \cite{CMP17,ChaMorNovPon19,CMNP19}, where a new notion of  distributional solutions was proposed.

	Let us focus on the definition of distributional solutions, referring to the nice review \cite{GigPoz22} for further information on viscosity solutions to \eqref{evol law} (we just note that the two notions are
	equivalent in the setting of this paper~\cite[Remark~6.1]{CMNP19}). 
	The exact definition of distributional solutions will be recalled in Definition~\ref{Defsol}, but  when $\phi$ is smooth it can be motivated as follows: It is known (see for instance \cite{Soner93} for the isotropic case) that $E(t)$ evolves according to \eqref{evol law} if and only if the  signed distance function $d(\cdot,t):=\sd^\po_{E(t)}$ to $\partial E(t)$  induced by the polar norm\footnote{defined by $\po(x)=\sup_{\phi(\nu)\le 1} \nu\cdot x$ and which satisfies
		$\phi(x)=\sup_{\po(x)\le 1} \nu\cdot x$.} $\po$, satisfies 
	\begin{align}
		\partial_t d\ge \div(\nabla \phi(\nabla d))\quad   \text{in }\{d>0\},\label{evol law 1}\\
		\partial_t d\le \div(\nabla \phi(\nabla d))\quad   \text{in }\{d<0\}\label{evol law 2}
	\end{align}
	in the viscosity sense. The idea of the new definition introduced in \cite{CMP17} is to reinterpret the equations above in the distributional sense. In particular, note that replacing $\nabla\phi(\nabla u)$ by 
	a vector field $z\in L^\infty(\{d>0\};\R^N)$ such that $z(x)\in \partial \phi(\nabla d)$ for a.e. $x$, where $\partial \phi$ denotes the subdifferential of $\phi$, the equations \eqref{evol law 1}, \eqref{evol law 2} make sense even when $\phi$ is crystalline. The corresponding notion of super- and sub-solutions bears a comparison principle, which yields uniqueness of the motion up to fattening. Existence is obtained either by a variant of the minimizing movements scheme of \cite{AlmTayWan,LucStu} in the spirit of  \cite{Chambolle04}, which consists in building a discrete-in-time evolution obtained by a recursive minimization procedure (see \cite{CMP17,CMNP19}), or	
	by approximation with smooth anisotropies~\cite{ChaMorNovPon19}.
	We observe that the convergence of such time discrete approaches to a motion characterized by~\eqref{evol law 1}-\eqref{evol law 2} in the \textit{viscosity sense} was shown in~\cite{KIshii14}, including
	in the 2D crystalline setting, while convergence in a distributional sense was established in~\cite{CasCha} in the convex case only.
	Briefly, given a time-step $h>0$ and an initial closed set $E_0=:E^{h,0}$, one defines $E^{h,k+1}=\{ u^{h,k+1}\le 0 \},$  where $u^{h,k+1}$ is defined as the minimizer of a so-called ``Rudin-Osher-Fatemi''~\cite{RuOsFa:92} problem:
	\begin{equation}\label{def u_k cont}
		u^{h,k+1}\in \text{argmin} \left\lbrace  \int_{\R^N} \phi(D u) + \frac 1{2h}\int_{\R^N} | u-\sd^\po_{E^{h,k}}|^2\right\rbrace.
	\end{equation}
	The idea of the present work
	is to combine this discretization in time with a simultaneous discretization in space for the particular class of purely crystalline anisotropies $\phi$ of the following form
	\begin{equation}\label{def phi}
		\phi(v)=\sum_{i\in \mathcal E} \beta(i)|i\cdot v|,
	\end{equation}
	where $\beta(i)>0$ and $\mathcal E\subseteq\Z^N\setminus\{0\}$ is a finite set of generators  such that  $\text{Span }\mathcal E=\R^N$.   {This kind of convex polytopes is known in the literature as \textit{rational zonotopes}. The class of rational zonotopes is dense in the class of symmetric convex sets if $N=2$, while for $N\ge 3$ it is nowhere dense. This fact is due to the strong symmetry properties of zonotopes, as every facet of a zonotope is itself a zonotope \cite{McM}. Note however that the Euclidean ball may be approximated by rational zonotopes in every dimension.}

	We now specify the discrete setting we are interested in, referring the reader to \cite{BraSol-Book} for a more thorough introduction to related topics. We consider  an $\e$-spaced square lattice $\ez$ and discrete functions $u : \ez\to \R$,  and denote $u_i:=u(i)$. We observe that we could also consider
	a general finite-dimensional Bravais lattice, at the expense of more tedious notation. 
	A natural discrete version of  total variation-like energies are those appearing in Ising systems, namely energies of the form 
	\begin{equation}
		TV_\beta^\e(v):=\e^{N-1}\sum_{i,j\in\ez} \beta(i/\e-j/\e)|v_i-v_j|,
	\end{equation}
	where $\beta$ is as in \eqref{def phi} and extended to 0 in $\Z^N\setminus \mathcal E$. 
	Under the hypotheses above on $\beta$, the functionals $TV_\beta^\e$  are shown to $\Gamma$-converge\footnote{Note that we do not need to assume that the lattice generated by $\{e_k\}_{k=1,\dots,m}$ is $\Z^N$, which is necessary to ensure the equi-coercivity of the discrete functionals.} as $\e\to 0$ to the total variation functional 
	\[
	TV_\phi(v)=\int_{\R^N}\phi(D v)
	\]
	where $\phi$ is as in \eqref{def phi}, see e.g.~\cite{ChaKre}. It is thus natural to define a minimizing movements scheme  based on $TV_\beta^\e$ which is the discrete counterpart of the minimizing procedure \eqref{def u_k cont}, as follows: given $E_0\subseteq \R^N,$ we define $E^{0}_{\e,h}=\{i\in\ez : (i+[0,\e)^N)\cap E_0\neq \emptyset\}$ and   for every $k\in\N$ we let $ u_{\e,h}^{k+1}$ be such that
	\begin{equation}\label{def u_k intro}
		u_{\e,h}^{k+1}\in \text{argmin} \left\lbrace  TV_\beta^\e(v)+\frac 1{2h}\sum_{i\in\ez}|v_i-(\sd^{k}_{\e,h})_i|^2 \ :\ v:\ez \to \R\right\rbrace,
	\end{equation}
	where $\sd^{k}_{\e,h}$ denotes a suitable signed $\po$-distance function to $E^k_{\e,h}$  defined  on $\ez$. (Actually, the energy in~\eqref{def u_k intro} is infinite and we rather
	consider the Euler-Lagrange equation of the problem.)
	Then, one sets $E^{k+1}_{\e,h}:=\{u^{k+1}_{\e,h}\le 0\}.$
	
	The idea is to study the asymptotic behaviour of the discrete evolutions  $E^{k}_{\e,h}$ as  both $\e,h\to 0$.
	A similar analysis has been performed
	in \cite{BraGelNov}, in the planar case, for $\phi=\|\cdot\|_1$ and $\sd^k_{\e,h}$ the continuous signed distance function from the discrete sets $E^{k}_{\e,h}$ restricted to the lattice $\ez$, see also \cite{MisiatsYip16,BraCicYip,BraSci,BraSol,MalNov,Sci} for further related results. With this choice, if $\e\gg h$ it is easy to see
	that the dissipation-like term in \eqref{def u_k intro}  
	\[
	\frac 1{2h}\sum_{i\in\ez}| v_i-(\sd^{k+1}_{\e,h})_i |^2
	\]
	forces the functions $u^k_{\e,h}$ to be constant as $k$ varies, therefore producing \textit{pinning} on the moving interfaces. 
	Moreover, when the two scales $\e,h$ are going to zero at the same speed it is shown in \cite{BraGelNov} that  a direct implementation of the standard scheme with the choice above for the distance, introduces
	a systematic error of order $\e=h$ at each step, which accumulates and produces a drift
	in the limiting evolution. As a result, low curvature shapes remain pinned, while sets with
	higher curvature evolve with a law which is a nonlinear modification of the crystalline curvature flow~\eqref{evol law}. Thus, the evolution law \eqref{evol law} can be approximated with the scheme of \cite{BraGelNov} only if $\e\ll h$. {In \cite{MisiatsYip16}, similar results are derived,
	still in dimension 2, for the isotropic (Euclidean) mean curvature flow.}
	
	We show in our main result, Theorem~\ref{themthm},
	that with a new appropriate definition of the distance $\sd^k_{\e,h}$, we can recover in the limit $\e,h\to 0$ the actual distributional solution to \eqref{evol law} for every initial set $E_0\subseteq \R^N$, for every purely crystalline anisotropy $\phi$ of the form~\eqref{def phi} with rational coefficients, in any dimension and irrespective of relative size of the space- and time-steps. In fact, the assumption of the rational character of $\beta$ can be removed in the regime $\e\leq O(h)$.  To the best of our knowledge this is the first general rigorous convergence result for a fully discrete scheme without restrictions on the dimension, on the initial sets and in which the spatial mesh is allowed to be of the same order or even coarser  than the time step.  
	
	Let us further comment on the  analysis carried out in \cite{BraGelNov} in the planar case (see also \cite{BraSol-Book} for many more references on the topic). One important change between these older results and ours is that we consider distributional solutions to the crystalline mean curvature flow \eqref{evol law}, instead of relying on the characterization of the motion via ODEs,  which dates back to \cite{AlmTay95,AngGu89}. The latter notion of solution is indeed suited only for planar evolutions, thus the limitation $N=2$ in the past works. 
	With the ODE definition and for $\phi=\|\cdot\|_1$, the authors of~\cite{BraGelNov} precisely prove the following results. If $\e\ll h$ then the limiting motion is consistent with \eqref{evol law}, while if $h\ll \e$ pinning happens for any nonempty initial data. As already mentioned, in the critical case $\e=h$, the limit planar motion is not driven by \eqref{evol law}, but instead by a slightly modified
	nonlinear crystalline mean curvature flow, and pinning may happen for some particular (low curvature) initial data. 
	This striking difference with our result may be (vaguely) justified by the following remark. While 
	in~\cite{BraGelNov}, the focus is on discrete sets, we rather evolve, in accordance with the
	definition of distributional solutions, the \textit{signed distance functions} to the  boundaries. In this way we can effectively achieve a sub-pixel precision in our approximation, as $u_{\e,h}$ and the signed distance function carry more information than the evolving level set $\{ u_{\e,h}(t)\le 0 \}$. Our new definition of the interpolated signed distance is detailed in Section~\ref{section MMS}.
	
	The consistency result in this paper validates the numerical experiments which we
	carry on in Section~\ref{sec:Num} to illustrate our results. These experiments are
	derived from previous experiments in~\cite{ChaDar}, which however were  using a
	different redistancing operation for which no consistency was proven. 
	Numerical schemes based on the variational approach~\cite{AlmTayWan,LucStu} have
	been introduced for crystal growth~\cite{Alm93}. Since then,
	there have been many attempts to implement implicit schemes based on this
	approach for isotropic and anisotropic curvature 
	flows in various settings~\cite{Chambolle04,EtoGigIsh,Obermanetal2011,Nor18,EtoGig}. { We are however not aware of a formal convergence proof for these schemes in the fully discrete setting, which is not}
	relying on the consistency of the spatial discretization with respect to the
	time-discrete scheme (hence assuming $\e\ll h$,
	{even if in practice these implementations seem very robust}). 
	
	Many other techniques have been considered to simulate crystalline flows after \cite{Tay92,Tay93}, see e.g. \cite{Gir95,GiraoKohn96,Dzi99} for the evolution of planar curves and  \cite{NovPao99,PaoPas00} for higher-dimensional algorithms.

	Let us conclude this introduction with two comments. The first one concerns the hypothesis that $\phi$ is purely crystalline. It seems quite technical, as it implies that the associated interaction function $\beta$ (in the sense of \eqref{def phi}) has finite range. While this is not necessary to carry {out} the existence part for the discrete minimizing movements scheme, it is essential for building a calibration
	which yields a bound on the speed of Wulff shapes, see Appendix~\ref{sec:proofscsv}. In practice, being the open  Wulff shape {$\W:=\{\po {\le}  1\}$} a finite Minkowski sum of (rational) segments (which is called a \textit{zonotope}), we can effectively handcraft a calibration along the directions identified by these segments. It is  a remarkable difference between this discrete setting and the continuous one, where instead the vector field $x/\po(x)$ in $\R^N$ is the right calibration \textit{for any} anisotropy $\phi$.
	
	The second one is on possible generalizations of the present analysis to more general evolution laws than \eqref{evol law}. The more general evolution law which is shown to admit a unique distributional solution is 
	\begin{equation}
		V(x,t)=\psi(\nu_{E(t)}(x))\left(-\kappa^\phi_{E(t)}(x)+f(x,t)\right), \qquad x\in\partial E(t), t\ge 0,
	\end{equation}
	where $\psi$ is a norm (usually referred to as the \textit{mobility}), and $f$ is a forcing term, see \cite{CMP17,CMNP19} . We expect most of the present analysis to be valid even if $\psi\neq \phi$, under suitable compatibility assumptions on $\psi$ (see \cite{CMP17,CMNP19} for details), and it should not be difficult to
	consider a driving force $f$ as long as it is Lipschitz in space and globally
	bounded, see~\cite{CMNP19} again.
	
	\bigskip
	
	The paper is organized as follows: in the next Section~\ref{subsec:cf}, we recall the
	definition of distributional crystalline curvature flows from~\cite{CMP17,CMNP19}. Then,
	we study the discrete ``Rudin-Osher-Fatemi'' problem and its Euler-Lagrange equation in
	Section~\ref{sec:ROF}. In Section~\ref{section MMS}, we introduce the discrete minimizing movement
	scheme, with our particular definition of the signed distance function. We study in detail
	the properties of these distances, then in Section~\ref{sec:Wulff} we analyse the particular
	case of an initial Wulff shape. In the continuous setting, it is well known that
	under the law~\eqref{evol law}, it decreases in a self-similar way with a speed proportional
	to the inverse of its radius. We show an estimate  bounding the decay of the discrete
	Wulff shapes, it relies on the delicate construction of a calibration $z$ for the Rudin-Osher-Fatemi problem with datum $\po$, detailed in Appendix~\ref{sec:proofscsv}.
	
	Our main result, which is that in the limit $\e,h\to 0$, 
	the motion defined in Section~\ref{section MMS} converges to a crystalline flow,
	is stated, and proved, in Section~\ref{sec:conv}. We implemented the discrete scheme
	in 2D and show some numerical simulations in Section~\ref{sec:Num}. Some technical
	results are collected in the Appendix.

	\section{Distributional crystalline curvature flows} 
\label{subsec:cf}
We recall the distributional formulation for the crystalline mean curvature
motion of sets evolving with normal velocity~\eqref{evol law}
introduced in \cite{CMP17} (see also \cite{CMNP19}). Here and in what follows $\phi$ is any norm, $\po$ denotes the polar (or dual) norm of $\phi$  and  given a closed set $F\subseteq\R^N$, $\dist^{\po}(\cdot, F)$ stands for the $\po$-distance function from $F$ defined by
\[
\dist^{\po}(x, F):=\min\{\po(x-y):\, y\in F\}\,.
\] 
Analogously, for any $E, F$ closed we set 
\[
\dist^{\po}(E, F):=\min\{\po(x-y):\, x\in E\,,\, y\in F\}\,.
\]
We recall that a sequence of closed sets $(E_k)_{k\ge 1}$ in $\R^N$ converges to a closed set $E$ in the {\em Kuratowski sense}:
if the following conditions are satisfied
\begin{itemize}
	\item[(i)] if $x_k\in E_k$ for each $k$, any limit point of $\{x_k\}$ belongs to $E$;
	\item[(ii)] for all $x\in E$ there exists a sequence $\{x_k\}$ such that $x_k\in E_k$ for each $k$ and $x_k\to x$.
\end{itemize}
We will write in this case:
\[
E_k\stackrel{\mathcal K}{\longrightarrow} E\,.
\]
One can easily verify  that $E_k\stackrel{\mathcal K}{\longrightarrow} E$ if and only if
(for any norm $\psi$) $\dist^{\psi}(\cdot, E_k)\to \dist^{\psi}(\cdot, E)$ locally uniformly in $\R^N$. Hence, by Ascoli-Arzel\`a Theorem we have that
any sequence of closed sets admits a  converging subsequence  in the Kuratowski sense
(possibly to $\emptyset$, when $\dist^{\psi}(\cdot, E_k)\to+\infty$).

\begin{definition}\label{Defsol}
	Let $E_0\subseteq\R^N$ be 
	a closed set. 
	Let $E$ be a closed set in $\R^N\times [0,+\infty)$ and
	for each $t\geq 0$ denote $E(t):=\{x\in \R^N\,:\, (x,t)\in E\}$. We
	say that $E$ is a {\em superflow} for \eqref{evol law} with
	initial datum $E_0$ if
	\begin{itemize}
		\item[(a)]$E(0)\subseteq {E}_0$;
		\item[(b)] 
		$E(s)\stackrel{\mathcal K}{\longrightarrow} E(t)$ as $s\nearrow t$ for all $t>0$;
		\item[(c)]
		If  ${E}(t)=\emptyset$ for some  $t\ge 0$, then $E(s)=\emptyset$ for all $s > t$.
		\item[(d)]  
		{Set $T^*:=\inf\{t>0\,:\, E(s)=\emptyset \text{ for $s\geq t$}\}$, 
			and 
			$$
			d(x,t):=\dist^{\po}(x, E(t)) \qquad \text{ for all } (x,t)\in \R^N\times (0,T^*)\setminus E.
			$$ }
		Then,
		\begin{equation}\label{eq:supersol}
			\partial_t d \ge \Div z 
		\end{equation}
		holds in the distributional sense in $\R^N\times (0,T^*)\setminus E$
		for a suitable $z\in L^\infty(\R^N\times (0,T^*))$ such that
		$z\in \partial\phi(\nabla d)$~a.e., $\Div z$ is a Radon measure in $\R^N\times (0,T^*)\setminus E$, and 
		{$(\Div z)^+\in L^\infty(\{(x,t)\in\R^N\times (0,T^*):\, d(x,t)\geq\delta\})$} for every $\delta\in (0,1)$.
	\end{itemize}
	
	We say that $A$, open set in $\R^N\times [0,+\infty)$, is
	a {\em subflow} for \eqref{evol law} with initial datum $E_0$ if $\R^N\times [0,+\infty)\setminus A$ is a superflow for \eqref{evol law} with initial datum $\R^N\setminus \textup{int}\big({E}_0\big)$.
	
	Finally, we say that $E$, closed set in $\R^N\times [0,+\infty)$,  is a {\em weak flow}  for \eqref{evol law} with initial datum $E_0$ if it is a superflow and if $\textup{int}(E)$\footnote{Here we are taking the interior with respect to $\R^N\times [0,+\infty)$} is a subflow, both with initial datum $E_0$.
\end{definition}

In \cite{CMP17} the following crucial inclusion principle between sub- and superflows is proven.

\begin{theorem}\label{th:compar}
	Let $E$ be a superflow  with initial datum $E_0$ and
	$F$ be a subflow with initial datum $F_0$ in the sense of Definition~\ref{Defsol}. Assume that
	$\dist^{\po}(E^0,\R^N\setminus{F^0})=:\Delta>0$. 
	Then, 
	$$\dist^{\po}\big(E(t),\R^N\setminus F(t)\big)\ge \Delta  \qquad \text{ for all } t\ge 0$$
	(with the convention that $\dist^{\po}(G, \emptyset)=\dist^{\po}(\emptyset, G)=+\infty$ for any $G$).
\end{theorem}

We  also recall the corresponding notion of sub- and supersolution to the level set flow associated with \eqref{evol law}. In what follows $\textup{UC}(\R^N)$ stands for the space of uniformly continuous functions on $\R^N$.

\begin{definition}[Level set subsolutions and supersolutions]\label{deflevelset1}
	Let  $u_0\in \textup{UC}(\R^N)$. A  lower semicontinuous   function $u:\R^N\times [0, +\infty)\to \R$ is called a  {\em  level set superflow}  for  \eqref{evol law}, with initial datum $u_0$,  if $u(\cdot, 0)\geq u_0$ and if for a.e. $\lambda \in \R$ the closed sublevel set $\{ u(\cdot,t)\leq \lambda\}$  is a superflow for  \eqref{evol law} in the sense of Definition~\ref{Defsol}, with initial datum $\{u_0\leq \lambda\}$.
	
	An  upper-semicontinuous   function $u:\R^N\times [0, +\infty)\to \R$ is called a  {\em  level set subflow} for \eqref{evol law}, with initial datum $u_0$, if~$-u$ is level set  superflow  in the previous sense, with initial datum $-u_0$.
	
	Finally,  a continuous   function $u:\R^N\times [0, +\infty)\to \R$ is called a {\em  level set flow} for  \eqref{evol law} if it is both a level set sub- and  superflow. 
\end{definition}
Using Theorem~\ref{th:compar}, it is not difficult to deduce the following parabolic comparison principle between level set sub- and superflows, which yields in particular the uniqueness of level set flows (in the sense of Definition~\ref{deflevelset1}), see \cite{CMNP19}. 
\begin{theorem}\label{th:lscomp}
	Let $u_0$, $v_0\in \textup{UC}(\R^N)$ and let $u$, $v$ be respectively a level set subflow starting from $u_0$ and a level set superflow starting from $v_0$. If $u_0\leq v_0$, then $u\leq v$.
\end{theorem}
We finally recall that in \cite{CMP17} (see also \cite{CMNP19}) the existence of level set flows is established by implementing a level-by-level minimizing movements scheme. This in turn yields
existence and uniqueness (up to fattening) for weak flows. This is made precise in the following  statement, see \cite[Corollary~4.6]{CMP17} and \cite[Theorem~4.8]{CMNP19}.
\begin{theorem}\label{th:phiregularlevelset} 
	Let  $u_0\in \textup{UC}(\R^N)$. Then the following holds: 
	
	{\rm (i)} There exists a unique level set flow $u$ in the sense of Definition~\ref{deflevelset1}   starting  $u_0$.
	
	{\rm (ii)}  For all $\lambda\in \R$ the sets $\{(x,t)\,:\, u(x, t) \le \lambda\}$ and $\{(x,t)\,:\, u(x, t) < \lambda\}$ are respectively the maximal superflow and minimal sublow with initial datum $\{u_0\leq \lambda\}$.
	
	{\rm (iii)}  For all but countably many $\lambda \in \R$, the fattening phenomenon does not occur; that is,  

	\begin{equation}
		\begin{array}{rcl}\label{eq:nonfattening}
			\{(x,t)\,:\, u(x, t) < \lambda\}& = &\mathrm{int\,}\big(\{(x,t)\,:\, u(x, t) \le \lambda\}\big)\,, \vspace{5pt} \\
			\mathrm{cl\,} \big(\{(x,t)\,:\, u(x, t) < \lambda\}\big) &=& \{(x,t)\,:\, u(x, t) \le \lambda\}\,,
		\end{array}
	\end{equation}
	where interior and closure are relative to space-time.
	
	For all such   $\lambda$, $\{(x,t)\,:\, u(x, t) \le \lambda\}$ is the unique  weak flow in the sense of Definition~\ref{Defsol}, starting from  $\{u_0\leq \lambda\}$.
\end{theorem}

The aim of this paper is to show that the convergence to the continuum level set flow holds true also when the Euler implicit time discretisation  is combined with a suitable spatial discretisation procedure.

\section{The discrete ``Rudin-Osher-Fatemi'' problem}\label{sec:ROF}

In this part, we describe our discrete setting, we then introduce and analyse
the discrete variant~\eqref{def u_k intro} of Problem~\eqref{def u_k cont}.

\subsection{Discrete functions spaces and operators} 
For $\e>0$, we define the function spaces $X_\e=\R^{\ez}$ and $Y_\e=\R^{\ez\times \ez}$.
Given a function $u\in X_\e$ and a discrete \virg{vector field} $ z\in Y_\e$, with a slight abuse of notation we will denote $u_i=u(i)$  and $z_{ij}=z(i,j)$, $i,j\in\ez$.
The discrete gradient $D_\e : X_\e\to Y_\e$ is defined, for $u\in X_\e$ as
\[ (D_\e u)_{ij}=\frac{u_i-u_j}{\e}.  \]
We denote its adjoint operator by $D^*_\e:Y_\e\to X_\e$, namely the operator such that, for $\eta \in X_\e$  compactly   supported and for $z\in Y_\e$, is defined as 
\[ \sum_{i} (D^*_\e z)_i\eta_i:= \sum_{ij} z_{ij} (D_\e\eta)_{ij}=\sum_{ij} z_{ij}{\frac{\eta_i-\eta_j}{\e}},   \]
where the indexes, here and throughout the paper, range over $\ez$ if not otherwise stated. In particular, taking $\eta=\chi_{\{i\}},$ one finds that
\begin{equation}\label{def D^*}
	(D_\e^* z)_i=\sum_{j} \frac{z_{ij}-z_{ji}}{\e},
\end{equation}
which can be seen as a discrete divergence operator.

\subsection{Discrete ROF problem}
In this section we consider the discrete anisotropic ROF problem associated with the discrete total variation functional. Without loss of generality, we consider $\e=1$ in this section, and
denote $X:=X_1$, $Y:=Y_1$ and $D:=D_1$. Given  a nonnegative $\beta\in X$, which will be called the \textit{interaction function}, satisfying
\begin{equation}\label{integr beta}
	\sum_{i\in\Z^N} \beta(i)=: c_\beta<+\infty,
\end{equation}
we set $\alpha_{ij}=\beta(i-j)$ and, for any $u\in X$ we define
\begin{equation}
	TV(u)=\sum_{i,j \in\Z^N} \alpha_{ij} |u_i-u_j| =\sum_{i,j} \alpha_{ij} |(Du)_{i,j}|.
\end{equation}
We also consider the discrete perimeter $\P$ defined for every $E\subseteq \Z^N$ as 
\[ \P(E):=TV(\chi^E)=\sum_{i,j\in\Z^N} \alpha_{ij} |\chi^E_i-\chi^E_j|.\]
We also consider a suitable localization of the perimeter: namely, for any set $A\subseteq \R^N$ we define 
\[ \P(E;A)=\sum_{i\in A\cap\Z^N \text{ or } j\in A\cap\Z^N} \alpha_{ij} |\chi^E_i-\chi^E_j|. \]
Note that the quantities above may well be infinite.

Then, given $g\in X$, we consider the following problem: find a pair $(u,z)\in X\times Y$ such that 
\begin{equation}\label{discrete EL}
	\begin{cases}
		D^* z + u = g\\
		z_{ij}(u_i-u_j)=\alpha_{ij}|u_i-u_j|, \quad |z_{ij}|\le \alpha_{ij}
		\quad \forall i,j\in\Z^N.
	\end{cases}
\end{equation}
The equation above is the Euler-Lagrange equation of the discrete ROF functional 
\begin{equation}\label{discrete ROF}
	ROF_g(v)=TV(v)+\frac 12 \sum_{i\in\Z^N}(v_i-g_i)^2.
\end{equation}
However, \eqref{discrete EL} makes sense also for those $g$ such that $ROF_g\equiv+\infty.$  
That~\eqref{discrete EL} is the first-order condition for optimality in~\eqref{discrete ROF} follows from standard convex
analysis: the idea is that since
\[
TV(v) = \sup\left\{ \langle z, Dv \rangle: |z_{i,j}|\le \alpha_{i,j} \forall (i,j)\right\},
\]
then the subgradients $\partial TV(v)$ of $TV$
at $v$ are precisely
given by the vectors $D^*z$ for the $z$'s which
realize the supremum in this expression. Then,
for $g$ with bounded support (such that there
is at least some $u$ with finite energy), \eqref{discrete EL}
expresses that $0\in\partial ROF_g(u)$, which by
definition is the condition for the minimality of $u$.

We will also consider the following  geometric minimization problem. Given $g\in X$, find
\begin{equation}\label{discr per pb}
	\min_{F\subseteq \Z^N} \P(F)+ \sum_{i\in\Z^N} \chi^F_i g_i.
\end{equation}

In order to deal with unbounded  sets, possibly with infinite perimeter, we will consider the following notion of global minimality with respect to compactly supported perturbations.

\begin{definition}\label{def loc minimizer}
	A set $E\subseteq \Z^N$ is a global minimizer for the problem \eqref{discr per pb} if for every $R>0$
	\begin{equation}\label{discr per pb loc}
		\P(E;B_R)+\sum_{|i|<R}\chi^E_i g_i\le  \P(F;B_R)+\sum_{|i|<R}\chi^F_i g_i
	\end{equation}
	for every $F\subseteq \Z^N$ such that $F\triangle E\subseteq B_R$. Here $B_R=\{x\in\R^N:|x|<R\}$
	is the open ball of radius $R$ centered in the origin.
\end{definition}

\begin{proposition}\label{geometric-comparison}
	Let $g,{g'}\in X$ such that $ g'-g\ge \delta >0$. Let  $E,E'$ be two global minimizers of problem  \eqref{discr per pb loc}, in the sense of Definition \ref{def loc minimizer}, corresponding to $g,g'$ respectively. Then, $E'\subseteq E.$
\end{proposition}
\begin{proof} 
	Let us denote in the following $\chi:=\chi^{E_s}, \chi':=\chi^{E'_s}$.  For a given $R>0$ we define the competitor sets $F=(E_s\setminus B_R) \cup ((E'_s\cup E_s)\cap B_R)$ and $F'=(E'_s \setminus B_R) \cup ((E'_s\cap E_s)\cap B_R)$. By minimality of $E_s, E'_s$ in $B_R$ one has
	\begin{align}
		\sum_{|i|<R \text{ or }|j|<R}\alpha_{ij} |\chi'_i-\chi'_j|+\sum_{|i|< R} g'_i(\chi'_i - \chi'_i\wedge \chi_i )&\le \sum_{\substack{|i|< R\\ |j|< R}}\alpha_{ij} |\chi'_i\wedge \chi_i-\chi'_j\wedge \chi_j| \label{comparison 1}\\
		&\ +  \sum_{\substack{|i| < R\\ |j| \ge R}}(\alpha_{ij}+\alpha_{ji}) |\chi'_i\wedge \chi_i-\chi'_j| \nonumber\\
		\sum_{|i|<R \text{ or }|j|<R}\alpha_{ij} |\chi_i-\chi_j|+\sum_{|i| <R} g_i(\chi_i - \chi'_i\vee \chi_i )&\le \sum_{\substack{|i| <R\\ |j| < R}}\alpha_{ij} |\chi'_i\vee \chi_i-\chi'_j\vee \chi_j|  \label{comparison 2}\\
		&\ +  \sum_{\substack{|i|< R\\ |j|\ge R}}(\alpha_{ij}+\alpha_{ji}) |\chi'_i\vee \chi_i-\chi_j| \nonumber.
	\end{align}
	Using the inequality\footnote{
		Indeed, if $a\ge b$ and $c\ge d$, this is an equality, while if $a>b$ and $c<d$, one deduces
		that $b-d < a-d < a-c$, $b-d < b-c < a-c$ so that there exists $t\in (0,1)$ with
		$a-d = t(b-d) + (1-t) (a-c)$, $b-c=(1-t)(b-d) + t(a-c)$: the conclusion follow by convexity of
		$|\cdot|$.
	} $|a\wedge b- c\wedge d|+ |a\vee b- c\vee d|\le |a-c|+|b-d|$ and summing together \eqref{comparison 1} and \eqref{comparison 2} we obtain
	\begin{equation}\label{comparison 3}
		\begin{split}
			&\sum _{\substack{|i| < R\\ |j| \ge R}} (\alpha_{ij}+\alpha_{ji}) \left( |\chi_i-\chi_j|+|\chi'_i-\chi'_j| \right) + 2\sum_{|i| < R} (g'_i-g_i)(\chi'_i - \chi_i)^+\\
			&\le  \sum_{\substack{|i|< R\\ |j|\ge R}} (\alpha_{ij} + \alpha_{ji}) \left( |\chi'_i\wedge \chi_i-\chi'_j| + |\chi'_i\vee \chi_i-\chi_j|  \right).
		\end{split}
	\end{equation}
	We then remark that $|\chi'_i\wedge \chi_i-\chi'_j|\le |\chi'_i\wedge \chi_i-\chi'_i|+ |\chi'_i-\chi'_j|=(\chi'_i-\chi_i)^+ + |\chi'_i-\chi'_j|$ and analogously $|\chi'_i\vee \chi_i-\chi_j|\le (\chi'_i-\chi_i)^+ + |\chi_i-\chi_j|$. Therefore, \eqref{comparison 3} entails 
	\begin{equation}\label{comparison 4}
		\sum_{|i| < R} (g'_i-g_i)(\chi'_i - \chi_i)^+ \le  \sum_{|i|< R} (\chi'_i - \chi_i)^+\sum_{|j|\ge R} (\alpha_{ij} + \alpha_{ji}).
	\end{equation}
	Fix now $R_\delta>0$ such that 
	\[ \sum_{|k|\ge R_\delta} \beta(k)\le \frac \delta4 \]
	and define $V_R:=\sum_{|i|< R} (\chi'_i - \chi_i)^+$. Assuming $R>R_\delta$, for every $\ell < R$ we use  \eqref{comparison 4} and $g+\delta\le g'$ to get   
	\begin{equation}
		\begin{split}\label{comparison 5}
			\delta V_R &\le \sum_{|i| < \ell} (\chi'_i - \chi_i)^+  \sum_{|j|\ge R} (\alpha_{ij} + \alpha_{ji}) + 2c_\beta \sum_{\ell \le |i|<R}  (\chi'_i - \chi_i)^+\\
			&\le 2\sum_{|i| < \ell} (\chi'_i - \chi_i)^+  \sum_{|k|\ge R-\ell} \beta(k) + 2c_\beta (V_R-V_\ell).
		\end{split}
	\end{equation}
	Therefore, choosing $\ell=R-R_\delta$ in \eqref{comparison 5} we obtain 
	\begin{equation}\label{iter 1}
		\frac \delta2 V_R\le 2c_\beta (V_R- V_{R-R_\delta}),
	\end{equation}
	which implies that for every $k,\ell\in \N$ it holds
	\begin{equation}\label{iter 2}
		V_{kR_\delta}\le \left( 1-\dfrac\delta{4c_\beta}   \right)^\ell V_{(k+\ell)R_\delta}.
	\end{equation}
	Letting $\ell\to +\infty$, since $V_{(k+\ell)R_\delta}=O(\ell^N),$ we infer that $V_{kR_\delta}=0 $ for every $k\in\N$.  In particular, this implies that $(\chi'-\chi)^+=0$ i.e. $\chi'\le \chi.$ 
\end{proof}

We will prove the following theorem.
\begin{theorem}\label{thm exist}
	Given $g\in X$ there exists a unique function $u^g\in X$ and there exists a discrete vector field $z\in Y$  such that $(u^g,z)$ is a solution of \eqref{discrete EL}. Moreover, the following comparison principle holds: if $g\le g'$ then $u^{g}\le u^{g'}$. Finally, for any $R>0$ and $s\in\R$ the sublevel set $E_s:=\{ i\in\Z^N : u_i^g\le s \}$ is a global minimizer (in the sense of Definition \ref{def loc minimizer}) for \eqref{discr per pb} with $g$ replaced by $g-s$.
\end{theorem}

\begin{proof}
	\textbf{Step 1.} (Existence)
	For every $n\in \N$ set $g^n:= g\chi^{B_n}$ and note that $g^n\in \ell^2(\Z^N)$.
	Therefore, by standard methods and by strict convexity the functional \eqref{discrete ROF}, with $g$ replaced by $g^n$ admits a unique minimizer $u^n$ and, as previously observed, the optimality condition is the existence of a discrete field $z^n$ such that $(u^n,z^n)$ solves \eqref{discrete EL} (with $g^n$ in place of $g$). Note that, for any $k\in\Z^N$, by equation \eqref{discrete EL} it holds
	\begin{equation}\label{equi L inf}
		|u_k^n|\le |g_k^n|+|(D^*z)_k|\le |g_k|+c_\beta\quad \text{for every }n\in\N,
	\end{equation}
	where the last inequality follows from the definition \eqref{def D^*} and  from $|z_{ij}|\le \alpha_{ij}$ and $|g^n|\le |g|$.
	Now, it is clear that we can extract a subsequence $n_k$ and find $(u,z)$ such that $u^{n_k}_i\to u_i$ and $z_{ij}^{n_k}\to z_{ij}$ as $k\to+\infty.$ Clearly we have that $|z_{ij}|\le \alpha_{ij}$ and $z_{ij}(u_i-u_j)=\alpha_{ij}|u_i-u_j|$ and it is immediate to check that $(u,z)$ satisfies equation \eqref{discrete EL}.\\
	\textbf{Step 2.} (Minimality of the sublevelsets) Let $R>0,s\in\R$ and let $F\subseteq \Z^N$ such that $E_s\triangle F\subseteq\joinrel\subseteq B_R.$ We first remark that $\alpha_{ij}|\chi^{E_s}_i- \chi^{E_s}_j|=-z_{ij}(\chi^{E_s}_i- \chi^{E_s}_j)$, which follows easily from the definition of $E_s$ and $z_{ij}(u_i-u_j)=\alpha_{ij}|u_i-u_j|.$
	
	We set $I_R:=\{ (i,j)\in \Z^N\times \Z^N : |i|<R \text{ or } |j|<R\}$ and compute
	\begin{equation}\label{calibra}
		\begin{split}
			&\P(F;B_R)-\P(E_s;B_R)=\sum_{(i,j)\in I_R} \alpha_{ij}|\chi^F_i- \chi^F_j| - \sum_{(i,j)\in I_R}  \alpha_{ij} |\chi^{E_s}_i- \chi^{E_s}_j| \\
			&\ge - \sum_{(i,j)\in I_R}  z_{ij}(\chi^F_i- \chi^F_j) +\sum_{(i,j)\in I_R}  z_{ij}(\chi^{E_s}_i- \chi^{E_s}_j)\\
			&= \sum_{(i,j)\in I_R}  z_{ij}(\chi^{E_s}_i-\chi^F_i- (\chi^{E_s}_j-\chi^F_j))\\
			&= \sum_{ij} z_{ij}(\chi^{E_s}_i-\chi^F_i- (\chi^{E_s}_j-\chi^F_j)),
		\end{split}
	\end{equation}
	where in the last equality we used the fact that $\chi^{E_s}_i=\chi^F_i$ if $|i|\ge R$. Noting that the function $\chi^{E_s}-\chi^F$ is compactly supported, we may use it as a test function for \eqref{discrete EL}. Therefore, from \eqref{calibra} we deduce
	\begin{equation*}
		\begin{split}
			&\P(F;B_R)-\P(E_s;B_R)\ge  \sum_{ij} z_{ij}(\chi^{E_s}_i-\chi^F_i- (\chi^{E_s}_j-\chi^F_j))\\
			&=\sum_i (\chi^{E_s}_i-\chi^F_i)(g_i-u_i)\ge \sum_{i\in E_s\setminus F} (g_i-s) -\sum_{i\in F\setminus E_s} (g_i-s),
		\end{split}
	\end{equation*}
	which shows the minimality of $E_s$.
	
	\textbf{Step 3.} (Comparison and uniqueness for \eqref{discrete EL})
	Assume $g\le g'$ and let $(u,z), (u',z')$ two corresponding solutions for \eqref{discrete EL}. Let $s>s'$ and recall that by Step 2 $\{ u'\le s'\}  $ and $\{ u\le s\}  $ are global minimizers for \eqref{discr per pb} according to Definition \ref{def loc minimizer}, with $g$ replaced by $g'-s'$ and $g-s$ respectively. Since $g'-s'-(g-s)\ge s-s'> 0$, from Proposition \ref{geometric-comparison} we obtain $\{ u'\le s'\}\subseteq \{ u\le s\}$. By the arbitrariness of $s,s'$ we conclude that $u\le u'.$
	
\end{proof}

\begin{remark}\label{rmk stup}
	We remark that, given $g\in X$ it clearly holds that $u^{-g}=-u^g$.
\end{remark}

\section{The minimizing movements scheme}\label{section MMS}
In this section we provide a combined spatial and time discretisation of the flow \eqref{evol law} for a particular class of 
norms $\phi$ and show the convergence of the scheme to the continuum flow. In what follows,  we consider $\{e_1,\dots,e_m\}\subseteq \Z^N$ a finite number of integer vectors spanning the whole $\R^N$, and set $\mathcal E=\{\pm e_k\}_{k=1}^m$. We let $\beta\in X$ be a non-negative function such that 
\[ \beta(-i)=\beta(i) \quad \text{ and }\quad \beta(i)>0\text{ if and only if  }i\in \mathcal E. \]
One can naturally associate an anisotropy $\phi$ with the function $\beta$ setting
\begin{equation}\label{phi special}
	\phi(v)=\sum_{i\in \mathcal E} \beta(i)|i\cdot v|= \sum_{k=1}^m 2\beta(e_k)|v\cdot e_k|.  
\end{equation}
Note that, in particular, it holds 
\begin{equation}\label{hp3}
	\#\{ k\in\Z^N : \beta(k)\neq 0\}<+\infty.
\end{equation}
We recall that the $\phi$-perimeter associated with \eqref{phi special} 
\[
P_\phi(E)=\int_{\partial^*E}\phi(\nu_E)\, d\mathcal{H}^{N-1}
\]
(defined for every $E\subseteq\R^N$ of finite perimeter) 
is the $\Gamma$-limit (in a suitable sense) as $\e\to 0$ of the following  scaled discrete perimeters 
\[
\P^\e(E):=\e^{N-1}\sum_{i,j\in\ez} \alpha^\e_{ij} |\chi^E_i-\chi^E_j|
= \e^N \sum_{i,j\in\ez} \alpha^\e_{i,j} |(D_\e \chi^E)_{i,j}|
\]
defined for all $E\subseteq \ez$, see for instance \cite{ChaBra23}.  Here  we have set
\begin{equation}\label{eq:alphah}
	\alpha^\e_{ij}:=\beta\left(\frac i{\e}-\frac j{\e}\right).
\end{equation}

Given $\phi$ a norm on $\R^N$ and a closed set $E\not\in\{\emptyset, \R^N\}$, we denote with $\sd^{\po}_E$ the  signed $\po-$distance function from $E$, which is defined as 
$$\sd^{\po}_E(x):=\min_{y\in E} \po(x-y)-\min_{y\notin E}\po(x-y). $$
We also set $\sd^{\po}_\emptyset\equiv +\infty$ and $\sd^{\po}_{\R^N}\equiv -\infty$.
We denote 
\begin{equation}\label{def c phi}
	C_\phi=\min_{i\in\Z^N\setminus \{0\}}\po(i)>0
\end{equation}
and define the $\phi$-Wulff shape $\W_R(x)$ of radius $R>0$ and center $x\in\R^N$ as $ \W_R(x)=\{ y\in\R^N : \po(x-y) {\le} R \}.$

\subsection{A discrete redistancing operator} In this section we introduce a discrete proxy for the signed distance function to a set, and study some of its properties.

Given $u\in X_\e$ we define the operators $d^{\e,\po}_\pm,\sd_\pm^{\e,\po},\sd^{\e,\po}: X_\e\to X_\e$ in the following way: letting $E=\{ i\in \ez : u_i\le 0 \},$ we first {set } 
\begin{equation}\label{def d,sd}
	\begin{split}
		(d^{\e,\po}_-(u))_i&=\sup_{j\in \{u\ge 0\}} \left\lbrace  u_j-\po(i-j) \right\rbrace,\\
		(\sd_-^{\e,\po}(u))_i&=\inf_{j\in  \{u\le 0\}}\left\lbrace  (d^{\e,\po}_-(u))_j+\po(i-j) \right\rbrace,\\
		(d^{\e,\po}_+(u))_i&=\inf_{j\in   \{u\le 0\}} \left\lbrace  u_j+\po(i-j) \right\rbrace,\\
		(\sd_+^{\e,\po}(u))_i&=\sup_{j\in  \{u\ge 0\}} \left\lbrace  (d^{\e,\po}_+(u))_j-\po(i-j) \right\rbrace,\\
		(\sd^{\e,\po}(u))_i&=\frac 12(\sd_+^{\e,\po}(u))_i+\frac 12 (\sd_-^{\e,\po}(u))_i.
	\end{split}
\end{equation}
Note that $d_+^{\e,\po}(u)=-d_-^{\e,\po}(-u)$ and $\sd_+^{\e,\po}(u)=-\sd_-^{\e,\po}(-u)$.

We will say that $f\in X_\e$ is $(L,\po)$-Lipschitz  if for all $i,j\in\ez$ it holds $|f_i-f_j|\le L\po(i-j)$. 

\begin{remark}\label{lip sd} We assume in what follows that $u$ is $(1,\po)$-Lipschitz.
	Then, concerning $ d^{\e,\po}_-,$ $ \sd^{\e,\po}_-$, we remark  that 
	\begin{equation}\label{def sd min}
		d^{\e,\po}_-(u)=\min\left\lbrace f\in X_\e : f \ge  u \text{ in } \{ u\ge 0 \},\  f \text{ is }  (1,\po)\text{-Lipschitz} \right\rbrace,
	\end{equation}
	and analogously
	\begin{equation}\label{def sd max}
		\sd^{\e,\po}_-(u)=\max\left\lbrace f\in X_\e : f \le d^{\e,\po}_-(u) \text{ in }  \{ u\le 0 \},\   f \text{ is }  (1,\po)\text{-Lipschitz} \right\rbrace.
	\end{equation}
	Correspondingly it holds 
	\begin{equation}\label{def sd+}
		\begin{split}
			&d^{\e,\po}_+(u)=\max\left\lbrace f\in X_\e : f \le u \text{ in }  \{u\le 0\},\   f \text{ is }  (1,\po)\text{-Lipschitz} \right\rbrace,\\
			&\sd^{\e,\po}_+(u)=\min\left\lbrace f\in X_\e : f \ge  d^{\e,\po}_+(u) \text{ in }  \{u\ge 0\},\  f \text{ is }  (1,\po)\text{-Lipschitz} \right\rbrace,
		\end{split}
	\end{equation}
	In particular, the functions $ d^{\e,\po}_\pm(u),\sd^{\e,\po}_\pm(u),\sd^{\e,\po}(u)$ are also $(1,\po)$-Lipschitz.
	Let us show \eqref{def sd min} the other identities being analogous. To this aim,  denote by $\hat d$ the function defined by the right-hand side of \eqref{def sd min}. Since $ d^{\e,\po}_-(u)$ is the pointwise supremum of   $(1,\po)$-Lipschitz  functions, we clearly have that  $ d^{\e,\po}_-(u)$ is itself  $(1,\po)$-Lipschitz. Moreover, testing with $j=i$ in the definition of $d^{\e,\po}_-(u)$, we get $ d^{\e,\po}_-(u)\geq u$ in  $\{ u\ge 0 \}$. Thus, we infer $\hat d\leq  d^{\e,\po}_-(u)$. For the opposite inequality, let $f$ be any functions as in the minimisation problem on the right-hand side of \eqref{def sd min}. Then for any $i\in \ez$ and $j\in \{u\geq 0\}$ we have
	\[
	f_i\geq f_j-\po(i-j)\geq u_j-\po(i-j)\,.
	\]
	By maximising with respect to $j\in \{u\geq 0\}$, we get $f\geq   d^{\e,\po}_-(u)$ and in turn, by the arbitrariness of $f$, $\hat d\geq   d^{\e,\po}_-(u)$, which concludes the proof of  \eqref{def sd min}
	
	Since the functions $ d^{\e,\po}_\pm(u),\sd^{\e,\po}_\pm(u),\sd^{\e,\po}(u)$ are  $(1,\po)$-Lipschitz, from \eqref{def sd min} it follows that
	\begin{equation}\label{d-=}
		d_-^{\e,\po}(u)\le u \  \text{in }\ez, \quad d_-^{\e,\po}(u)= u \  \text{in }\{u\geq 0\},
	\end{equation}
	while  \eqref{def sd max} implies that
	\begin{equation}\label{sd-=}
		\sd^{\e,\po}_-(u)\ge d_-^{\e,\po}(u) \ \text{in }\ez, \quad \sd^{\e,\po}_-(u)= d_-^{\e,\po}(u) \  \text{in }\{u\leq 0\}.
	\end{equation}
	Reasoning in the same way, we see that  
	\begin{equation}\label{sd+=}
		\begin{split}
			d_+^{\e,\po}(u)\ge u \  \text{in }\ez&, \quad d_+^{\e,\po}(u)= u \  \text{in }\{u\leq 0\}, \\
			\sd^{\e,\po}_+(u)\le d_+^{\e,\po}(u)  \  \text{in }\ez&, \quad \sd^{\e,\po}_+(u)= d_+^{\e,\po}(u) \  \text{in } \{u\geq 0\}.
		\end{split}
	\end{equation}
	In particular we conclude
	\begin{equation}\label{sd=}
		\sd^{\e,\po}(u)\geq  u \ \text{in }\{u\geq 0\}, \quad \sd^{\e,\po}(u)\leq  u \  \text{in }\{u\leq 0\}.
	\end{equation}
	{Note that \eqref{sd=} implies $\{\sd^{\e,\po}_+(u)\ge 0\} \supseteq \{u \ge 0\}$, and    \eqref{sd+=} yields $\{\sd^{\e,\po}_+(u)< 0\} \supseteq \{u < 0\}$, thus $\{\sd^{\e,\po}_+(u)\ge 0\} = \{u \ge 0\}$ (and analogously for $\sd^{\e,\po}_-$). Similarly, one shows that $\{\sd^{\e,\po}_\pm(u)\le 0\} = \{u \le 0\}$.}
	In particular, if the level set $0$ of $u$
	is ``fat'', then this is preserved by these discrete ``signed distance functions''.
	Further properties of these discrete signed distance functions are presented in Lemma~\ref{lem:sdcompar} below  and in Remark \ref{rmk:sd}
	
	Moreover, it follows directly from the definition of $d^{\e,\po}_\pm(u),\sd^{\e,\po}_\pm(u)$ that the function $\sd^{\e,\po}(u)$ is invariant under integer translations, meaning that for any $i,\tau\in\ez$ it follows 
	\begin{equation}\label{transl invar sd}
		\left(\sd^{\e,\po}(u(\cdot+\tau))\right)_i=\left(\sd^{\e,\po}(u)\right)_{i+\tau}.
	\end{equation}
\end{remark}

{We now show that the redistancing operator $\sd^\e(u)$ is indeed a discrete approximation of the signed distance function to the $0$-sublevel set of the function $u$.}

{Given a set $E\subseteq \ez,$ we will denote with $\widehat E\subseteq \R^N$ the closed set defined by }
\[\widehat E:=E+[0,\e]^N. \]

\begin{lemma}\label{lemma sd sd hat}
	Given a $(1,\po)$-Lipschitz function $u\in X_\e$, {it holds}
	\begin{equation}\label{comp sd sd hat}
		\sup_{\ez\setminus E}|\sd^{\e,\po}_\pm(u)-\sd_{\widehat E}^{\po}|\le c_\phi \e,
	\end{equation}
	for a suitable positive constant $c_\phi$,  where  $E=\{i\in \ez : u_i\le 0 \}$. Moreover, 
	\begin{equation}\label{comp sd sd hat bis}
		\sd^{\e,\po}_\pm(u)\ge \sd_{\widehat E}^{\po}- c_\phi \e \quad \text{in }\ez.
	\end{equation}
\end{lemma}

\begin{proof}
	In this proof we let $c_\phi$ denote a positive constant which depends on $\phi$ and that may change from line to line and also within the same line.

	We start introducing a slightly modified definition of the discrete signed distance $\sd^{\e,\po}(u)$. Namely, setting
	\begin{equation}\label{def boundary+}
		\begin{split}
			&\pa^+_\e E:=\{ i\in \ez\setminus E \ : \  \exists j\in E \text{ with }\|i-j\|_\infty=\e\} \\
			&\pa^-_\e E:=\{ i\in E \ :\  \exists j\in \ez\setminus E \text{ with }\|i-j\|_\infty=\e\}
		\end{split},
	\end{equation}
	we define
	\begin{equation}
		\tilde d_i=  \begin{cases}
			\inf \left\lbrace   u_j + \po(i-j)\ :\ j\in \pa^-_\e E \right\rbrace,\quad &\text{for }i\in \ez\setminus E\\
			\sup \left\lbrace   u_j - \po(i-j)\ :\ j\in \pa^+_\e E  \right\rbrace & \text{for }i\in E
		\end{cases} .
	\end{equation}
	We start by showing that
	\begin{equation}\label{eq comp 3}
		\begin{split}
			\sd_\pm^{\e,\po}(u)\ge \tilde d\quad &\text{in }E,\\
			\sd_\pm^{\e,\po}(u)\le \tilde d\quad &\text{in }\ez\setminus E.
		\end{split}
	\end{equation}
	Indeed, we note that for every $i\in E$ we have 
	\[ (\sd^{\e,\po}_-(u))_i=(d_-^{\e,\po}(u))_i=\sup_{j\in \{u\ge 0\}}   \{ u_j-\po(i-j) \}\ge \sup_{j\in\pa^+_\e E}   \{ u_j-\po(i-j) \} = \tilde d_i. \]
	On the other hand, recalling that $d_-^{\e,\po}(u)\le u$ in $ E$,  for every $i\in \ez\setminus E$ we see
	\[ (\sd^{\e,\po}_-(u))_i=\inf_{j\in  \{u\le 0\}}   \{ (d_-^{\e,\po}(u))_j+\po(i-j) \}\le \inf_{j\in  \pa_\e^- E}  \{ u_j+\po(i-j) \} = \tilde d_i.\]
	Reasoning analogously  we show the same inequalities between $\sd^{\e,\po}_+$ and $\tilde d$ and thus prove  \eqref{eq comp 3}.

	Next, we prove
	\begin{equation}\label{tilde d and sd}
		\sup_{\ez}|\tilde d-\sd_{\widehat E}^{\po}|\le c_\phi \e. 
	\end{equation}
	Recall that by definition \eqref{def boundary+}, since $u\le 0$ in $E$ and $u> 0$ in $\ez\setminus E$ and since $u$ is $(1,\po)$-Lipschitz, it holds 
	\[ |u_j|\le c_\phi \e \quad \text{for } j\in\pa_\e^\pm E.\]
	Then, for every $i\in \ez\setminus E$ we have
	\begin{equation}\label{eq comp 1}
		\tilde d_i=\inf_{j\in \pa^-_\e E}  \{   u_j + \po(i-j)\}\ge \inf_{j\in \pa^-_\e E}  \po(i-j)-c_\phi \e \ge \sd_{\widehat E}^{\po}(i) -c_\phi \e. 
	\end{equation}
	On the other hand, by definition of $ \sd_{\widehat E}^{\po}$ there exists $x\in\pa \widehat E$ such that $\sd_{\widehat E}^{\po}(i)=\po(i-x)$. Let $k\in\ez$ be the closest point from $x$ in $\pa^-_\e E$. We have
	\begin{equation}\label{eq comp 2}
		\begin{split}
			\sd_{\widehat E}^{\po}(i)&=\po(i-x)\ge \po(i-k)-c_\phi \e\\
			&\ge \po(i-k)+ u_{k}-c_\phi \e\ge \tilde d_i- c_\phi \e.
		\end{split}
	\end{equation}
	Finally, equation \eqref{eq comp 1} and \eqref{eq comp 2} imply \eqref{tilde d and sd} outside $E$. The other case is analogous.

	We now finally prove \eqref{comp sd sd hat} outside $E$. From \eqref{eq comp 3} and \eqref{tilde d and sd} it holds \begin{align*}
		d_-^{\e,\po}(u)=\sd^{\e,\po}_-(u)\ge \tilde d\ge \sd_{\widehat E}^{\po}- c_\phi \e \quad \text{ in }E.
	\end{align*}
	In particular, $\sd_{\widehat E}^{\po}- c_\phi \e$ is an admissible competitor in \eqref{def sd max}, thus $ \sd^{\e,\po}_-(u)\ge  \sd_{\widehat E}^{\po}- c_\phi \e$  in $\ez$.
	On the other hand, in $\ez\setminus E$ it holds \eqref{eq comp 3}, thus we conclude \eqref{comp sd sd hat} for  $\sd^{\e,\po}_-(u)$. Concerning $\sd^{\e,\po}_+(u)$, we note that  by Remark~\ref{lip sd} and the equation above it holds
	\[ u\ge \sd^{\e,\po}_-(u)\ge \sd^\po_{\widehat E}-c_\phi \e\quad \text{in }E. \]
	The function $\sd^\po_{\widehat E}-c_\phi \e$ is therefore admissible in \eqref{def sd+}, thus by maximality 
	\[ d^{\e,\po}_+(u)\ge \sd^\po_{\widehat E}-c_\phi \e. \]
	Since $\sd^{\e,\po}_+(u)=d^{\e,\po}_+(u)$ in $\ez\setminus E$  we conclude \eqref{comp sd sd hat}, taking also into account  again \eqref{eq comp 3} and \eqref{tilde d and sd}. Finally, \eqref{comp sd sd hat bis} follows by combining \eqref{comp sd sd hat}, \eqref{eq comp 3} and \eqref{tilde d and sd}.
\end{proof}

We conclude the section with some further properties of the operator $\sd^{\e,\po}$.

\begin{lemma}\label{lem:sdcompar}
	Given $u\in X_\e$ and $(1,\po)$-Lipschitz, it holds 
	\begin{equation}\label{simm sd}
		\sd^{\e,\po}(-u)=-\sd^{\e,\po}(u).
	\end{equation}
	Furthermore, if $u_1,u_2\in X_\e$ are $(1,\po)$-Lipschitz and $u_1\le u_2$ then 
	\begin{equation}\label{sd monotone}
		\sd^{\e,\po}(u_1)\le \sd^{\e,\po}(u_2).
	\end{equation}
	Finally, for any $s>0$ and $u\in X_\e$ and $(1,\po)$-Lipschitz, it holds
	\begin{equation}\label{sd level sets}
		\sd^{\e,\po}(u-s)\le \sd^{\e,\po}(u)-s. 
	\end{equation}
\end{lemma}

\begin{proof}
	For every $i\in \ez$ it  holds
	\begin{equation*}
		(d^{\e,\po}_-(-u))_i=\max_{ j\in \{(-u)\ge 0\}} \left\lbrace  -u_j-\po(i-j) \right\rbrace=-\min_{ j\in \{u\le 0\}} \left\lbrace  u_j+\po(i-j) \right\rbrace=-(d^{\e,\po}_+(u))_i.
	\end{equation*}
	In turn, 
	\begin{equation*}
		\begin{split}
			(\sd^{\e,\po}_-(-u))_i&=\min_{ j\in \{(-u)\le 0\}} \left\lbrace  (d_-^{\e,\po}(-u))_j+\po(i-j) \right\rbrace\\
			&=- \max_{ j\in \{u\ge 0\}} \left\lbrace  (d_+^{\e,\po}(u))_j-\po(i-j) \right\rbrace=-(\sd^{\e,\po}_+(u))_i.
		\end{split}
	\end{equation*}
	Reasoning in the same way for $d^{\e,\po}_+,\sd^{\e,\po}_+$  we arrive at
	\begin{equation}\label{sd pm}
		\sd^{\e,\po}_\pm(-u)= -\sd^{\e,\po}_\mp(u)
	\end{equation}
	and thus $ \sd^{\e,\po}(-u)= -\sd^{\e,\po}(u)$. The monotonicity property \eqref{sd monotone} follows easily from Definition~\eqref{def d,sd}.
	The proofs of the other results also follow  from Definition~\eqref{def d,sd}, we present only the one concerning \eqref{sd level sets}. Fix $s>0$ and $u\in X_\e$ be a  $(1,\po)$-Lipschitz function. By definition of $ d^{\e,\po}_-(u)$ we have
	\[ (d^{\e,\po}_-(u))_i=\sup_{ j\in \{u\ge 0\}} \left\lbrace  u_j-\po(i-j) \right\rbrace \ge s+ \sup_{ j\in \{u\ge s\}} \left\lbrace  (u_j-s)-\po(i-j) \right\rbrace =(d^{\e,\po}_-(u-s))_i+s. \]
	Analogously
	\begin{align*}
		(\sd^{\e,\po}_-(u))_i&=\inf_{ j\in \{u\le 0\}} \left\lbrace  (d_-^{\e,\po}(u))_j+\po(i-j) \right\rbrace \\
		&\ge s+\inf_{ j\in \{u\le s\}} \left\lbrace  (d_-^{\e,\po}(u-s))_j+\po(i-j) \right\rbrace=s+(\sd^{\e,\po}_+(u-s))_i.
	\end{align*}
	Since the proofs for $d^{\e,\po}_+(u),\sd^{\e,\po}_+(u)$ are analogous, we conclude.
\end{proof}

\subsection{The discrete scheme}
We now describe our minimizing movements scheme,  discretized in both time and space. 
A particularity of our scheme is that in practice, it evolves the distance function
to a set rather than the set itself. In particular, at the discrete level, it may
depend on the initialization (even if in the limit the flow is geometric and only
depends on the initial set).

Recalling~\eqref{eq:alphah}, we  rescale equation \eqref{discrete EL} on the lattice $\ez$ in the following way.  We recall that $X_\e=\R^{\ez}$ and $Y_\e=\R^{\ez\times \ez}$. 
Given $g\in X_\e$ {and a time step $h>0$,} the problem \eqref{discrete EL} now becomes  to  find $  (u,z)\in X_\e\times Y_\e$   satisfying 
\begin{equation}\label{discrete EL resc}
	\begin{cases}
		h D^*_\e z +  u=g \quad \text{on }\ez\\
		z_{ij}(u_i-u_j)=\alpha^\e_{ i  j }|u_i-u_j|, \quad |z_{ij}|\le \alpha^\e_{ i  j },
	\end{cases}
\end{equation}
where $D^*_\e z$ is defined  in \eqref{def D^*}.
{For ease of notation we assume $\e=\e(h)$, with $\e\to 0$ as $h\to 0$ and  we will specify the dependence on $h$ only.}

Let $E_0\subseteq \R^N$ be a closed set. We define $E^{h,0}:=\{ i\in \ez : (i+[0,\e)^N)\cap E_0\neq \emptyset \}.$ We note that
\begin{equation}
	\widehat E^{h,0}\to E_0, \quad E^{h,0} \to E_0 
\end{equation}
as $h\to 0$ in the Kuratowski sense,  where with a slight abuse of notation we write $\widehat E^{h,0}$ to denote the set $\widehat{ E^{h,0}}=E^{h,0}+[0,\e]^N$.

Given a closed set $E_0\subseteq \R^N$ with $E_0\notin \{\emptyset,\R^N\}$, we
consider $u^{h,0}$ a $(1,\po)$-Lipschitz function on $\ez$ which is negative
inside $E^{h,0}$ and positive outside. For instance, we set
\[ u^{h,0}:=\frac12 C_\phi \e(1-\chi_{E^{h,0}})-\frac12 C_\phi \e \chi_{E^{h,0}}, \]
where $C_\phi$ is defined in \eqref{def c phi}, so that $u^{h,0}$ is $(1,\po)$-Lipschitz. Let us set $(z^{h,0})_{ij}=0$ for all $i,j\in\ez$. Then, as long as $E^{h,k}\notin \{\emptyset,\R^N\}$, we can iteratively define  $u^{h,k+1},z^{h,k+1}$ for $k\in\N$ by solving \eqref{discrete EL resc} with $g=\sd^{\e,\po}(u^{h,k})$; i.e., 
\begin{equation}\label{discrete EL resc algo}
	\begin{cases}
		h D^*_\e z^{h,k+1} +  u^{h,k+1}=\sd^{\e,\po}(u^{h,k}) \quad \text{on }\ez\\
		z^{h,k+1}_{ij}(u^{h,k+1}_i-u^{h,k+1}_j)=\alpha^\e_{ i  j }|u^{h,k+1}_i-u^{h,k+1}_j|, \quad |z^{h,k+1}_{ij}|\le \alpha^\e_{ i  j }.
	\end{cases}
\end{equation}
{We recall that the redistancing operator $\sd^{\e,\po}$ has been introduced in the previous section.} We then set
\[ E^{h,k+1}=\{ i\in \ez : u^{h,k+1}_i\le 0 \}. \]
If either $E^{h,k}=\emptyset$ or $E^{h,k}=\R^N$, we define $E^{h,k+1}=E^{h,k}.$ We denote by $T^*_h$ the first discrete time $hk$ such that $E^{h,k}=\emptyset$, if any; otherwise we let $T^*_h=+\infty$. Analogously, we set ${T'}^{*}_h$ first discrete time $hk$ such that $E^{h,k}=\R^N$, if any; otherwise we let $T'^{*}_h=+\infty$. 

For ease of notation we will set 
\begin{equation}\label{def algo}
	\begin{split}
		&E^h(t):=E^{h,[t/h]}\subseteq \ez\\
		&d^{h}(t):=\sd^{\e,\po}(u^{h,[t/h]})\in X_\e\\
		&u^h(t):=u^{h,[t/h]}\in X_\e\\
		&z^h(t):=z^{h,[t/h]}\in Y_\e\\
		&\widehat d^{h}(\cdot,t):=\sd^{\po}_{\widehat E^h(t)}\in \text{Lip}(\R^N),
	\end{split}
\end{equation}
where again, with a slight abuse of notation, $\widehat E^h(t)$ stands for $\widehat {E^h(t)}$. Note that in the definition of $\widehat d^{h}(\cdot,t)$ we are possibly using the convention  $\sd^{\po}_{\emptyset}\equiv+\infty$ and $\sd^{\po}_{\R^N}\equiv-\infty$. Note also that  $z^h(t)$ is well defined only for $0\leq t<\min\{T^*_h, T'^{*}_h \}$; however, if needed, we can set $z^h(t)=0$ for $t\ge \min\{T^*_h, T'^{*}_h \}$.   

\begin{remark}\label{lip u}
	If $u$ is the solution of \eqref{discrete EL resc} with  $(L,\po)$-Lipschitz datum $g$, by standard arguments, based on  the comparison principle and  translation invariance, one can show that $u$ satisfies the same Lipschitz bound of $g$. Indeed, given $j\in \ez$, the function $u(\cdot-j)\pm L\po(j)$ solves \eqref{discrete EL resc} with datum $g(\cdot-j)\pm L\po(j)$. By comparison one concludes as $ g(\cdot-j)-L \po(j) \le g(\cdot)\le g(\cdot-j) + L \po(j) .$
\end{remark}

\begin{lemma} 
	Let $u^{h}, E^{h}, d^{h}$ be defined as in \eqref{def algo}. Then, {for every $t\ge 0$, $d^{h}(t)$} is $(1,\po)$-Lipschitz and satisfies 
	\begin{equation}\label{u over d }
		\begin{cases}
			u^{h}(t)\le d^{h}(t) \quad &\text{in } \ez\setminus E^{h}(t)\\
			u^{h}(t)\ge d^{h}(t) &\text{in }E^{h}(t).
		\end{cases}
	\end{equation}
\end{lemma}

\begin{proof}
	It follows from Remarks \ref{lip sd} and \ref{lip u}.
\end{proof}

\begin{remark}\label{rmk evol compl} (Evolution of the complement)
	Let $E^h(t),u^h(t)$ be as in \eqref{def algo}. We note that, if $F_0\subseteq \R^N$ is a closed set such that $F^{h,0}=\ez\setminus E^{h,0}$, then the discrete evolution starting from $F_0$ coincides with  $ \{ u^h(t)\ge 0 \}$  for every $t\ge 0$. Indeed, denoting $v^h$ the discrete evolution starting from $F_0$, it holds by definition $v^{h,0}=-u^{h,0}$, thus recalling   \eqref{simm sd} we have
	\[  \sd^{\e,\po}(v^{h,0})=-\sd^{\e,\po}(u^{h,0}) \]
	and, by uniqueness for \eqref{discrete EL resc} it follows that $v^h(h)=-u^h(h)$. Then we can iterate to conclude.
\end{remark}

\begin{remark}[Comparison principle]\label{rmk comparison}
	Let $E_0, F_0$ be closed sets in $\R^N$ such that $E^{h,0}\subseteq F^{h,0}$ (note that this condition is satisfied if $E_0\subseteq F_0$).  Let $E^h(t),F^h(t)$ be the corresponding discrete evolutions and let $u^h(t),v^h(t)$ be the associated functions as in \eqref{def algo}. Then, for every $t\ge 0$ it holds $E^h(t)\subseteq F^h(t)$. This follows easily by iteration  from the monotonicity property \eqref{sd monotone} and from the comparison principle for \eqref{discrete EL resc}. One in fact could also consider the \virg{open} discrete evolution given by $ \mathring E^h(t):=\{ u^h(t)<0 \}$ and $ \mathring F^h(t):=\{ v^h(t)<0 \}$. Then, by the same argument one also have that $\mathring E^h(t)\subseteq\mathring F^h(t).$
\end{remark}

\begin{remark}[Avoidance principle]\label{rmk avoidance}
	Let $E_0,F_0\subseteq\R^N$ be closed sets such that $ E^{h,0}\cap  F^{h,0}=\emptyset$ (which is, for example, implied by $\textup{dist}(E_0,F_0)> c_\phi \e$ for a suitable $c_\phi>0$). Let $E^h,u^h$ and $\mathring F^h(t),v^h$ be  the  closed and open discrete evolutions starting from $E_0,F_0$ respectively (where the open discrete evolution has been defined in Remark \ref{rmk comparison}). Then, 
	\[\mathring F^h(t)\subseteq \ez\setminus E^h(t).\]
	Indeed, $F^{h,0}\subseteq \ez\setminus E^{h,0}$ implies that $-u^{h,0}\le v^{h,0}$ and thus by \eqref{simm sd} and \eqref{sd monotone}
	\[ -\sd^{\e,\po}(u^{h,0})=\sd^{\e,\po}(-u^{h,0})\le \sd^{\e,\po}(v^{h,0}). \]
	By the comparison principle for \eqref{discrete EL resc} and iterating one sees that $-u^h(t)\le v^h(t)$ for all $t\ge 0,$ which implies 
	\[  \mathring F^h(t)=\{ v^h(t)<0 \}\subseteq \{ u^h(t)>0 \}=\ez\setminus E^h(t). \]

\end{remark}

\begin{remark}\label{rmk:sd}
	We conclude this section by observing that we could have made different choices of the distance function, without affecting the final convergence result. In  definition \eqref{def d,sd} we could have set 
	\begin{equation}\label{def sd weird}
		\begin{split}
		    ( d^{< }(u))_i&=\inf_{j\in  \{u< 0\}} \left\lbrace  u_j+\po(i-j) \right\rbrace,\\
			( \sd^{<}(u))_i&=\sup_{j\in \{u\ge 0\}} \left\lbrace  ( d^{<}(u))_j-\po(i-j) \right\rbrace,\\
			( d^{\le }(u))_i&=\inf_{j\in  \{u\le 0\}} \left\lbrace  u_j+\po(i-j) \right\rbrace,\\
			( \sd^{\le}(u))_i&=\sup_{j\in \{u> 0\}} \left\lbrace  ( d^{<}(u))_j-\po(i-j) \right\rbrace.
		\end{split}
	\end{equation}
    One can see that $\sd^\le(u) $ mimics the signed distance function to the boundary of $\{u \le 0\}$ while  $\sd^<(u) $ mimics the signed distance function to the boundary of $\{u < 0\}$. Defining the algorithm as in \eqref{discrete EL resc algo} but with $\sd^<,\sd^\le$ replacing $\sd^{\e,\po}$, adapting our proof one can conclude the same convergence result. Let us further comment on the relation between $\sd^{\e,\po},\sd^\le,\sd^<$. One can prove that for any $(1,\po)$-Lipschitz function $u\in X_\e$, then
	\begin{equation}\label{three sd}
		\sd^\le(u) \le \sd^{\e,\po}_-(u)\le \sd^{\e,\po}_+(u)\le \sd^<(u).
	\end{equation}
	Thus, between the many possible choices we could have performed in \eqref{def d,sd}, it turns out that $\sd^<$ is the \virg{maximal} one, while $\sd^\le$ is the \virg{minimal}. Indeed, let us show that $\sd^{\e,\po}_-(u)\le \sd^{\e,\po}_+(u).$ By definition \eqref{def d,sd} and \eqref{d-=}, \eqref{sd+=} for every $i\in \{u\ge 0\}$ it holds 
	\[ (\sd^{\e,\po}_-(u))_i=\inf_{j\in \{u\le 0\}} \left\lbrace  (d^{\e,\po}_-(u))_j+\po(i-j) \right\rbrace {\le}  \inf_{j\in \{u\le 0\}} \left\lbrace  u_j+\po(i-j) \right\rbrace = (\sd^{\e,\po}_+(u))_i. \]
	Reasoning analogously, for every $i\in\{u\le 0\}$ it holds 
	\[ (\sd^{\e,\po}_+(u))_i=\sup_{j\in \{u\ge 0\}} \left\lbrace  (d^{\e,\po}_+(u))_j-\po(i-j) \right\rbrace {\ge}  \sup_{j\in \{u\ge 0\}} \left\lbrace  u_j-\po(i-j) \right\rbrace = (\sd^{\e,\po}_-(u))_i. \]
	Furthermore,  for any two $(1,\po)$-Lipschitz functions $u,u'\in X_\e$, if $u\le u'-s$ for $s>0$ then \[ \sd^<(u)\le \sd^\le(u')-s.  \]
	In particular, this implies that for any  $(1,\po)$-Lipschitz function $u\in X_\e$ and $s'>s$ then 
	\[ \sd^{\e,\po}(u-s)\le \sd^{\e,\po}(u-s') +s'-s. \]
	Fix  $u_0\in X_\e$ is a $(1,\po)$-Lipschitz function. Using the properties above and standard arguments, one can see that for all but countably many $s\in \R$ the discrete evolutions starting from  $\{u_0\le s\}$ and corresponding to the three  possible choices of distances in \eqref{three sd} coincide.
\end{remark}

\subsection{Discrete evolution of Wulff shapes}\label{sec:Wulff}
In this section we provide some control on the evolution speed of discrete Wulff shapes. The first result estimates the solution of~\eqref{discrete EL resc} for the distance to the Wulff shape.
\begin{lemma}\label{speriamo che sia vero}
	There exists a constant $C=C(\phi)>0$ with the following         property. If $u$ is the solution of \eqref{discrete EL resc} with $g=\po$, then  $u\le \phi^h$, where $\phi^h\in X_\e$ is defined as
	\begin{equation}\label{phi h}
		\phi^h_i:=\begin{cases}
			\phi^\circ(i)+\frac {Ch}{\phi^\circ(i)}    \quad &\text{if }\phi^\circ(i)\ge C (\sqrt h \vee \e)\\
			C (\sqrt h \vee \e) + \frac{Ch}{\sqrt h \vee\e }  &\text{otherwise}.
		\end{cases}
	\end{equation}
\end{lemma}
The proof of Lemma~\ref{speriamo che sia vero}, based on the construction of a calibration, is postponed to Appendix~\ref{sec:proofscsv}. We now prove a useful lemma used to estimate the redistancing step in our algorithm for functions of the form of \eqref{phi h}.
\begin{lemma}\label{lemma sd po-R}
	Let $R\ge \delta>0$ and set 
	\[
	u:=(\po-R)\vee(\delta/2-R).
	\]
	Then,
	for $\e$  small enough depending on $\delta$ it holds
	\begin{equation}
		\sd^{\e,\po}(u)\le \po-R +\hat c \e\quad \text{in }\ez,
	\end{equation}
	for a suitable positive constant $\hat c$, depending on $\phi$. 	Furthermore, if we assume \eqref{rational zonotope}, it holds 
	\begin{equation}\label{sd po-R}
		\sd^{\e,\po}(u)\le \po-R \quad \text{in }\ez.
	\end{equation}
\end{lemma}
\begin{proof}
    By \eqref{three sd}, it is sufficient to prove the claim for $\sd^{\e,\po}_+$.  We start showing that $d^{\e,\po}_+(u)=u$, 
     noting that by \eqref{sd+=} it suffices to prove $d^{\e,\po}_+(u)\le u$ in $\{ u\ge0 \}=\{ \po\ge R \}.$  	Assuming \eqref{rational zonotope}, given $i\in \{u\ge 0\}$  we note that $\po(i)\ge R$ thus by Lemma \ref{lemma neigh} there exists $j\in \W_R\setminus \W_{R-2\e\ell_1}$ satisfying
	\begin{equation*}
		\po(j)+\po(i-j)=\po(i).
	\end{equation*} 
	Taking $\e=\e(\delta)$ we can ensure that $R-2\e\ell_1 \ge \delta/2$, so that $j\in (\W_R\setminus \W_{\delta/2})\cap \ez$.  By definition {\eqref{def d,sd}} and the equation above we conclude that
	\[ d^{\e,\po}_+(u)\le u_{ j} + \po(i-j)= \po(j)-R +  \po(i-j)= \po(i)-R,\]
	hence we have shown that $d^{\e,\po}_+(u)=u$.
	Finally, from the definition \eqref{def d,sd} and since $d^{\e,\po}_+(u)=u=\po-R$ on $\{u\ge 0\}$, we conclude by the triangular inequality that $\sd^{\e,\po}_+(u)\le \po-R$. All in all, we have obtained \eqref{sd po-R}.
	
	If instead \eqref{rational zonotope} does not hold, using the first part of Lemma \ref{lemma neigh} and reasoning as above, one concludes that
	\[
	\sd^{\e,\po}_+(u)\le \po-R+\hat c\e,
	\]
	for a positive constant $\hat c$, and then the conclusion follows.
\end{proof}

Combining the two results above we can provide a bound on the evolution speed of Wulff shapes in the algorithm \eqref{discrete EL resc algo}.

\begin{proposition}\label{prop:evol Wulff}
    Assume either $\e\le O(h)$ or that \eqref{rational zonotope} holds. For every $\delta>0$ there exist $\e_0,h_0,c_0$ positive constants depending on $\delta$ with the following property. If $R\ge \delta$, $\e\le \e_0$ and $ h\le h_0$, then the discrete evolution of $\W_R$ defined in  \eqref{discrete EL resc algo}, denoted $\W^h(t)$,  satisfies 
    \begin{equation}\label{est2 speed wulff}
        \W^h(t) \supseteq (\W_{R-c_0 (t+\e)})\cap \ez,
    \end{equation}  
	as long as  $R- c_0(t+\e) \ge \delta/2$.
\end{proposition}

\begin{proof}
     Let   $\mathring \W^{h}(t)$ be the open discrete  evolution (see Remark \ref{rmk comparison})  starting from the closure of $\W_R$, for some $R>0$  and let $v^h(t)$ be the associated function as in {the third equation in} \eqref{def algo}.  Using the definition of $v^{h,0}$, \eqref{sd-=} and the first definition in  \eqref{def d,sd}, it is easy to see that
	\begin{equation}\label{sd=-R:1}
		(\sd^{\e,\po}_-(v^{h,0}))_0=(d^{\e,\po}_-(v^{h,0}))_0\le  -R+c_\phi\e.
	\end{equation}
	On the other hand, consider $i\in \{v^{h,0}\ge 0\}$ and let $ x'\in \partial \W_R$ be such that
	\[ \po(i-x')= \po(i)-\po(x')=\po(i)-R. \]
	Since there exists $j'\in  \{v^{h,0}\le 0\}$ such that $\po( j' - x')\le c_\phi \e$, then by triangular inequality 
	\[  \po(i- j')\le \po(i)-R+c_\phi \e. \]
	Thus, using again definition \eqref{def d,sd}, we get
	\begin{equation*}
		(d^{\e,\po}_+(v^{h,0}))_i\le\inf_{j\in  \{v^{h,0}\le 0\}}\po(i-j)\le \po(i)-R+c_\phi\e,
	\end{equation*}
	which implies
	\begin{equation}\label{sd=-R:2}
		(\sd^{\e,\po}_+(v^{h,0}))_0\le\sup_{j\in  \{v^{h,0}\ge 0\}}  (d^{h,0}_+(v^{h,0}))_j-\po(j)\le-R+c_\phi \e.
	\end{equation}
	Therefore, since $\sd^{\e,\po}(v^{h,0})$ is a $(1,\po)$-Lipschitz function, from \eqref{sd=-R:1}, \eqref{sd=-R:2} we get that 
	\[ \sd^{\e,\po}(v^{h,0})\le \po -R+c_\phi\e \quad \text{in }\ez. \]
	By comparison and Lemma \ref{speriamo che sia vero} we obtain 
	\begin{equation}\label{start}
		v^h(h)\le \phi^h -R +c_\phi \e,
	\end{equation}
	where $\phi^{h}\in X_\e$ is defined  in \eqref{phi h}. Considering $R\ge \delta$ and $h=h(\delta), \e=\e(\delta)$ small enough, the equation above implies that
	\begin{equation}\label{start2}
		v^h(h)\le \left(\po - R +  c_0h +c_\phi\e \right)\vee \left(\frac{\delta}{2}-R\right)
	\end{equation}
	where $c_0=4C/\delta$, with $C$ the same as in \eqref{phi h}. 	Assume first \eqref{rational zonotope}. From Lemma \ref{lemma sd po-R}, with $R$ replaced by $R-c_0 h-c_\phi \e$, we get 
	\begin{equation}\label{stima 1}
		\sd^{\e,\po}(v^h(h))\le \po -R + c_0 h+ c_\phi \e,
	\end{equation}
	therefore by comparison and Lemma \ref{speriamo che sia vero} we get 
	\[
	v^h(2h)\le \phi^h-R + c_0 h+c_\phi\e,
	\]
	which, reasoning as above, implies for $\e(\delta),h(\delta)$ small
	\[
	v^h(2h)\le (\po -R + 2 c_0 h+c_\phi\e)\vee \left( \frac \delta2 -R\right).
	\]
	Hence, we can iterate the argument to conclude that 
	\begin{equation}\label{est speed vh}
		v^h(t)\le (\po-R +  c_0 t + c_\phi\e)\vee \left( \frac \delta2 -R\right),
	\end{equation}
	as long as  $R -  c_0 t -  c_\phi \e \ge \delta/2$ and $\e,h$ are sufficiently small.
	In particular, this implies \eqref{est2 speed wulff} (possibly changing the value of $c_0$).
	
	If instead \eqref{rational zonotope} does not hold and $\e \le O(h)$,  we obtain \eqref{start}, \eqref{start2}  in the same way. Then, using the first part of Lemma \ref{lemma sd po-R} we get 
	\begin{equation}
		\sd^{\e,\po}(v^h(h))\le \po -R + c_0 h+\hat c\e+ c_\phi \e,
	\end{equation}
	then iterating we get
	\[
	v^h(kh)\le (\po-R + kc_0 h+ k\hat c \e +c_\phi\e)\vee\left( \frac\delta 2 -R \right),
	\]
	hence, recalling that $\e\le O(h)$ we conclude \eqref{est speed vh} and \eqref{est2 speed wulff}, 
	as long as  $R -  c_0 t-c_\phi\e\ge \delta/2$, with $\e,h$ sufficiently small and possibly changing the value of $c_0$.	
\end{proof}

As a corollary of the previous result, we deduce an estimate of the evolution of the distance function $\widehat d^h$ at distance from the evolving boundary, which we show next.

\begin{corollary}\label{lm:bound dist}
	Let $E_0\subseteq \R^N$ be a closed set and consider the discrete evolution defined in \eqref{def algo}. Assume either that $\e\le O(h)$ or that \eqref{rational zonotope} holds.
	Then, for every $\delta>0$ there exist $c_0=c_0(\delta)>0$, $h_0=h_0(\delta)>0$ and $\e_0=\e_0(\delta)$ such that the following holds.  If  $\widehat d^h(x,t)\ge  \delta,$ then for $s\ge t$, 
	\begin{equation}\label{bound dist+}
		\widehat d^h(x,s) \ge \widehat d^h(x,t) - c_0 (s-t+\e+h)
	\end{equation}
	provided $0<h\le h_0$, $0<\e<\e_0$ and as long as $\widehat d^h(x,t) - c_0 (s-t+\e+h) \ge \delta/2$. Similarly, if  $\widehat d^h(x,t)\le -\delta,$ then for $s\ge t$, 
	\begin{equation}\label{bound dist-}
		\widehat d^h(x,s) \le \widehat d^h(x,t) + c_0 (s-t+\e+h)
	\end{equation}
	provided $0<h\le h_0$ and as long as  $ \widehat d^h(x,t) + c_0 (s-t+\e+h)\le -\delta/2$.
\end{corollary}

\begin{proof}
	As usual, in this proof we denote by $c_\phi$ a positive constant depending on $\phi$ whose value may change from line to line and also within the same line.
	
	Assume $\widehat d^h(x,t)\ge \delta$. Without loss of generality we may assume $t\in [0, T^*_h)$ so that $\widehat d^h(x,t )$ is finite. Denote by $x_\e\in \ez$ such that $x\in x_\e+[0, \e)^N$. Note that there exists a constant $c_\phi>0$ such that, setting 
	$R:= \widehat d^h(x,t)- c_\phi \e$, one has $(\W_R(x_\e))^{h,0}\cap E^h(t)=\emptyset$   and $R>\delta/2$ 
	(if $\e,h$ are sufficiently small, depending on $\delta$).  By the avoidance principle stated in Remark \ref{rmk avoidance}, we deduce that the open discrete evolution of $\W_R(x_\e)$, which we denote by  $F(\tau$), lies outside $E^h([\frac th]h+\tau)$ for all $\tau\ge 0.$ By Proposition \ref{prop:evol Wulff} we deduce
	\begin{equation}\label{incl wulff}
		F(\tau)\supseteq \W_{R-c_0(\tau+\e)}(x_\e)\cap \ez,
	\end{equation}
	provided that $R-c_0 (\tau+\e)\ge \delta/2$. Note that in particular 
	\[ (\W_{R-c_0 (\tau+h+\e)}(x_{{\e}})\cap \ez) \subseteq (\ez\setminus E^h(t+\tau)), \]
	as long as $R-c_0 (\tau+h+\e) \ge \delta/2$. In turn, we get 
	\begin{equation}
		\widehat  d^h(x_\e,t+\tau)\ge R-c_0(\tau+h+\e),
	\end{equation}
	provided $R-c_0 (\tau+h+\e) \ge \delta/2$ (for a  possibly larger value of $c_0$). Recalling the definition of $R$ and $x_\e$ and possibly increasing the value of $c_0$, we  infer 
	\begin{equation}
		\widehat  d^h(x,t+\tau)\ge  \widehat  d^h(x,t)-c_0(\tau+h+\e)
	\end{equation}
	as long as $\widehat  d^h(x,t)-c_0(\tau+h+\e)\ge \delta$. The case $\widehat d^h(x,t) \le -\delta$ is analogous.
\end{proof}

\section{Convergence of the scheme}\label{sec:conv}
We now are ready to study the convergence of the scheme as $\e\to 0, h \to 0$. Recall that we assumed that $\e=\e(h)$ goes to 0 as $h\to 0$. In this section we assume that either $\e\le O(h)$ or that \eqref{rational zonotope} holds.
Let  $E^h(\cdot)$ be the discrete evolution defined in \eqref{def algo} and recall that $\widehat E^h(\cdot)=E^h(\cdot)+[0,\e]^N$. We introduce the closed space-time tubes 
\begin{equation}\label{tubidiscreti}
	\overline E^h:=\textup{cl}\big(\{(x,t)\in \R^N\times[0,+\infty):\, x\in \widehat E^h(t)\}\big)\,
\end{equation}
where the closure is in space-time. 
Then,  there exist $A, E$ open and closed (respectively) subsets of $ \R^N\times[0,+\infty)$, with $A\subseteq  E$, and a  subsequence $h_k\to 0$ such that
\[
\overline E^{h_k}\stackrel{\mathcal K}{\longrightarrow} E\qquad\text{and}\qquad 
\R^N\times[0,+\infty)\setminus\textup{int}\big( \overline {E}^{h_k}\big)\stackrel{\mathcal K}{\longrightarrow}  \R^N\times[0,+\infty)\setminus A,
\]
where interior, and Kuratowski convergence are  meant in space-time.
Let $E(t)$ and $A(t)$ be the $t$-time slice of $E$ and $A$, respectively.. 

Note that if $E(t)=\emptyset$ for some $t\ge 0$,
then \eqref{bound dist+} implies $E(s)=\emptyset$ for all $s\ge t$
so that we can define, as in Definition~\ref{Defsol}, the extinction
time $T^*$ of $E$. In the same fashion  one can define the extinction  time ${T'}^*$ of
$ \R^N\times[0,+\infty)\setminus A$ (notice that at least one between  $T^*$ and ${T'}^*$ is $+\infty$).
Possibly extracting a further (not relabelled) subsequence and arguing exactly as in \cite[Proof of Proposition 4.4]{CMP17} (and relying on the bounds \eqref{bound dist+} and \eqref{bound dist-}), one can in fact show the following result.

\begin{proposition}\label{prop:E}
	There exists a countable set $\mathcal N\subseteq (0, +\infty)$ such that
	${\widehat d^{h_k}}(\cdot, t)^+\to \dist^{\po}(\cdot, E(t))$   and  $\widehat d^{h_k}(\cdot, t)^-\to \dist^{\po}(\cdot,\R^N\setminus A(t))$ locally uniformly  for all $t\in (0, +\infty) \setminus \mathcal{N}$.
	Moreover, $ E$ and $ \R^N\times[0,+\infty)\setminus A$ satisfy the continuity properties (b) and (c) of Definition~\ref{Defsol}. In addition,  if $T^*>0$, then  $\{{\widehat d^{h_k}}\}$ is locally uniformly bounded in $\R^N\times (0, T^*)\setminus E$ and  analogously   $\{{\widehat d^{h_k}}\}$  is locally uniformly bounded in $\R^N\times (0, {T'}^*)\cap A$ if ${T'}^*>0$.
	Finally, $E(0)= E_0$ and  $A(0)=\textup{int}(E_0)$. 
\end{proposition}

\begin{theorem}\label{themthm}
	The set $E$ is a superflow in the sense of Definition~\ref{Defsol}
	with initial datum $E_0$,
	while $A$ is a subflow with initial datum ${E}_0$.
\end{theorem}

The proof of this result follows the main lines of the proof of~\cite[Theorem 4.5 ]{CMP17}.
One important difference with respect to the local, continuous setting is that the variable
$z^{h_k}$ is defined on the \textit{edges} $(i,j)$ between the vertices ${i,j}\in \ez$ and
it is therefore unclear how to pass to the limit in this variable to obtain the limiting
vector field $z(x,t)$. 
In order to do so, we associate with the discrete vector field $z_{ij}^h(t)\in Y_\e$ a vector field $\mathbf{z}^h(\cdot, t)$ in $\R^N$  defined as follows:
\begin{equation}
	\mathbf{z}^h(x, t):=\frac 1\e \sum_{j\in\ez} z_{ij}^{h}(t)(i-j),
\end{equation}
where $i\in\ez$ is such that  $x\in i+[0,\e)^N$. Recall that we can take $z_{ij}^{h}(t)$ and thus $ \mathbf{z}^h(\cdot, t)$ identically zero for $t\ge \min\{T^*_h, T'^{*}_h \}$. First, we show
the following:
\begin{lemma}\label{lemma phi circ}
	The vector field $\mathbf{z}^h$ satisfies
	\begin{equation}
		\po(\mathbf{z}^h)\le 1.
	\end{equation}
\end{lemma}

\begin{proof}
	Take $v\neq 0$ in $\R^N$. Recalling that $\phi(v)=\sum_{\ell\in \Z^N} \beta(\ell) |v\cdot \ell|$, one has for any $x\in\R^N$ and $i\in\ez $ such that $x\in i+[0,\e)^N$
	\begin{equation}\label{eq:defzcontinu}
		\mathbf{z}^h(x,t)\cdot v=\frac 1\e \sum_{j\in\ez} z_{ij}^{h}(t)(i-j)\cdot v= \sum_{\ell \in \Z^N} z_{i,i+ \e \ell }^{h}(t) \ell \cdot v\le \phi(v),
	\end{equation}
	where we used that $|z_{i,i + \e \ell}^{h}(t)|\le \beta(\ell).$
\end{proof}

Hence, being globally bounded, this vector field is weakly-$*$ compact in  $L^\infty(\R^N\times (0,T);\R^N)$ for any $T>0$. The following lemma establishes a relationship between the divergence of its limits and the limits of the discrete divergences of $z^h$.

\begin{lemma}\label{lemma d star to div}
	Assume that $\mathbf{z}^{h_k}\overset{*}{\rightharpoonup} z$ in $L^\infty(\R^N\times (0,T);\R^N)$ along a   subsequence $h_k\to 0$. Then,  for every $\varphi\in C^\infty(\R^N\times (0,T))$ and $\eta\in C^\infty_c(\R^N\times (0,T))$ it holds 
	\[ \lim_{k\to \infty} \left(  \e_k^N  \int \sum_{i,j\in \kz}   z^{h_k}_{ij}(t)\eta(i,t)\,\frac{\varphi(i,t) -\varphi(j,t)}{\e_k} \ud t\right) = \iint   \eta\,z\cdot\nabla\varphi\ud x\ud t. \]
\end{lemma}

\begin{proof}
	Let $\varphi\in C^\infty(\R^N{\times(0,T)})$ and $\eta\in C^\infty_c(\R^N{\times(0,T)})$ and denote $S(t)=\textup{supp}(\eta(t))$ and $Q_k:=[0,{\e_k})^N$. We have
	\begin{equation}\label{end align1}
		\e_k^{N} \sum_{i,j\in\kz}  z^{h_k}_{ij}(t)\eta (i,t)\, \frac{ \varphi(i,t)-\varphi(j,t)}{\e_k} = \e_k^N\sum_{i,j\in\kz} \frac {z^{h_k}_{ij}(t)}{\e_k} \eta(i,t)\, \nabla \varphi(x_{ij})\cdot (i-j),
	\end{equation}
	where $x_{ij}$ belongs to the segment between $i$ and $j$. 
	Furthermore we have
	\begin{align}
		&\left\lvert  \e_k^{N} \sum_{i,j\in\kz} z^{h_k}_{ij}(t)\eta(i,t)\, \frac{ \varphi(i,t)-\varphi(j,t)}{\e_k} -\sum_{i,j\in \kz}\frac {z^{h_k}_{ij}(t)}{\e_k}  \int_{i+Q_k}\eta \nabla  \varphi\cdot (i-j) \ud x \right\rvert \nonumber\\
		&\le \sum_{i,j\in \kz} \frac {\alpha^{\e_k}_{ij}}{\e_k} |\eta(i,t)|  \int_{i+Q_k} \left\lvert \left(  \nabla \varphi(x_{ij},t)-  \nabla \varphi(x,t)\right) \cdot (i-j)\right\rvert  \ud x   +    O(\e_k^{N})\label{eq align2}\\
		&\le 2\|\eta\|_{\infty}\sum_{i\in S(t)\cap \kz}  \sum_{j\in \kz}   \frac {\alpha_{ij}^{\e_k}}{\e_k}   \int_{i+Q_k} \left\lvert \left(  \nabla \varphi(x_{ij},t)-  \nabla \varphi(x,t)\right) \cdot (i-j)\right\rvert  \ud x+ O(\e_k^{N})\nonumber\\
		&\le c \e_k^N \sum_{i\in  S{(t)}\cap \kz}  \sum_{j\in \kz}   \frac {\alpha_{ij}^{\e_k}}{\e_k}|i-j|^2 +  O(\e_k^{N})\label{ineq lip varphi}\\
		&=c \e_k^{N+1}\sum_{\substack{i\in \Z^N\\
				\e_k  i\in S(t)}}\sum_{j\in \Z^N}  \alpha_{ij} |i-j|^2+ O(\e_k^{N})\nonumber\\
		&\le c \e_k^{N+1} \left( \sum_{\ell\in\Z^N} \beta(\ell) |\ell|^2 \right) (\# S(t)\cap \kz)+ O(\e_k^{N})\nonumber\\
		&\le c \e_k \sum_{\ell\in\Z^N} \beta(\ell) |\ell|^2+ O(\e_k^{N})\label{end align}
	\end{align}
	where in \eqref{eq align2} we used the Lipschitz property of $\eta$ and  \eqref{hp3}, while in \eqref{ineq lip varphi} we used the Lipschitz property of $\nabla \varphi$ and  $|x_{ij}-x|\le (1+\sqrt N)|i-j|$ for $i\neq j$ and $x\in i+Q_k$, and finally in \eqref{end align} we used that $ \#(S(t)\cap \ez)=O(\e_k^{-N}),$ which holds locally uniformly in time. Moreover, note that the the estimate provided above is uniform as $t$ varies in compact subsets of $(0,T)$.  Recalling \eqref{hp3}, we conclude integrating in time and sending $k\to \infty$.
\end{proof}

At this point, we may proceed with the proof of Theorem~\ref{themthm}.
\begin{proof}[Proof of Theorem~\ref{themthm}]
	As usual, in this proof we denote by $c_\phi$ a positive constant depending on $\phi$ whose value may change from line to line and also within the same line.
	
	We only show that $E$ is a superflow, as the subflow property of $A$ can be proven analogously.
	Points (a), (b) and (c) of Definition~\ref{Defsol} follow from Proposition~\ref{prop:E}.
	We are left with showing (d). Without loss of generality we may assume $T^*>0$ (which follows from Corollary \ref{lm:bound dist} if the initial set is not trivial). Note also that by Proposition~\ref{prop:E} we have $\liminf_kT^*_{h_k}\geq T^*$. \\
	\noindent\textbf{Step 1:} (Proof of \eqref{eq:supersol}). 
	For $(x,t)\in \R^N\times (0, T^*)\setminus E$ we set $d(x,t):=\dist^{\po}(\cdot, E(t))$. By Lemma~\ref{lemma sd sd hat} and Proposition~\ref{prop:E} we have 
	\begin{equation}\label{conv1}
		\sup_{\kz\cap K}|d^{h_k}(t)-d(\cdot, t)|\to 0\text{ as $k\to\infty$ for $t\in (0, T^*)\setminus \mathcal N$ and for any compact $K\subseteq \R^N\setminus E(t)$.}
	\end{equation}
	Moreover, $d^{h_k}$ and $d$ are locally uniformly bounded in $\R^N\times (0, T^*)\setminus E$.
	Set  $\z^{h_k}(\cdot, t):=0$ for $t>T^*_{h_k}$ if $T^*_{h_k}<T^*$. Extracting a further subsequence, if needed, and recalling Lemma~\ref{lemma phi circ},
	we may assume
	that $\z^{h_k}$ converges  weakly-$*$ in $L^\infty(\R^N\times (0,T^*); \R^N)$  to some vector-field $z$
	satisfying
	\begin{equation}\label{subdiff 0}
		\po(z)\le 1
	\end{equation}
	almost everywhere.
	Recall that by~\eqref{u over d } we have $u^{h_k}(t) \le d^{h_k}(t)$ in $\kz\setminus E^{h_k}(t)$; i.e., in the region where $d^{h_k}(t)$ is nonnegative. Combining with \eqref{discrete EL resc algo} (and recalling \eqref{def algo}) we infer that for $t<T^*_{h_k}$ it holds
	\begin{equation}\label{eq:ineqdisc}
		-D^*_{\e_k}z^{h_k} (t+h_k) \le \frac{d^{h_k}(t+h_k) -d^{h_k}(t)}{h_k} \qquad \text{ in  } \kz\setminus E^{h_k}(t).
	\end{equation}
	Consider a nonnegative test function $\varphi\in C_c^\infty((\R^N\times (0,T^*))\setminus E)$. If $k$ is large enough, then the distance of the support of $\varphi$  from $\overline E^{h_k}$ is bounded away from zero. In particular,  $d^{h_k}$ is finite  and positive on $\mathrm{\textup{supp}\,}\varphi$. We deduce from~\eqref{eq:ineqdisc} that 
	\begin{align}
		&\e_k^N\int \sum_{i\in\kz} \varphi(i,t)\left(\frac{d^{h_k}_i(t+{h_k})-d^{h_k}_i(t)}{h_k}+(D^*_{\e_k}z^{h_k}(t+h_k))_i \right)
		\ud t \nonumber \\ 
		&=-\e_k^N\int  \sum_{i\in\kz}  \frac{\varphi(i,t)-\varphi(i,t-h_k)}{h_k}d^{h_k}_i(t) \ud t   +   \e_k^N\int \sum_{i,j\in\kz} \frac{z^{h_k}_{ij}(t+h_k)-z^{h_k}_{ji}(t+h_k)}{h_k} \varphi(i,t)
		\ud t\nonumber\\ 
		&=-\e_k^N\int  \sum_{i\in\kz}  \frac{\varphi(i,t)-\varphi(i,t-h_k)}{h_k}d^{h_k}_i(t) \ud t    +   \e_k^N\int  \sum_{i,j\in\kz} z^{h_k}_{ij}(t+h_k) \frac{\varphi(i,t)-\varphi(j,t)}{h_k} 
		\ud t \ge 0.\label{eq end align}
	\end{align}
	It is easy to check that the first integral in \eqref{eq end align} converges to $-\iint d\,\partial_t \varphi \ud x\ud t  $ as $k\to \infty$ thanks to \eqref{conv1} and since $d^{h_k}, d$ are uniformly bounded. Recalling that $\z^{h_k}$  converges weakly-$*$ in $L^\infty(\R^N\times (0,T^*))$ to $z$, we use Lemma \ref{lemma d star to div} to conclude that the second integral in \eqref{eq end align} converges to $\iint z\cdot \nabla \varphi\ud x\ud t$. We thus conclude  \eqref{eq:supersol}.
	\\
	\noindent\textbf{Step 2:} (Convergence of $u^{h_k}$ to $d$).
	Firstly, we establish an upper bound for $ -D^*_{\e_k} z_{h_k}$ away from $E^{h_k}$. 
	We start by noting that definition \eqref{def d,sd}   implies
	\begin{equation}\label{bound Lip sd}
		\sd^{\e,\po}(u)\le \frac12\left( (d_-^{\e,\po}(u))_j+u_\ell+\po(\cdot-j) + \po(\cdot-\ell)  \right)\quad \text{in }\ez\setminus \{u\le 0\},
	\end{equation}
	for every $(1,\po)$-Lipschitz function $u\in X_\e$ and $j,\ell\in \{ u\le 0\}.$ Therefore, specifying the inequality above for $u^{h_k}(t)$,  by the comparison principle  and Lemma~\ref{speriamo che sia vero} we conclude
	\begin{equation}
		u^{h_k}_i(t+h_k)\le \frac12\left( \phi^{h_k}_{i-j}+\phi^{h_k}_{i-\ell} +(d_-^{\e,\po}(u^{h_k}(t))_j+u^{h_k}_\ell(t)  \right), \quad \forall i\in\kz\setminus E^{h_k}(t),
	\end{equation}
	where $j,\ell\in E^{h_k}(t)$.  If $\widehat d^{h_k}(i,t)\geq R>0$, recalling the definition of $\phi^h$, we get 
	\begin{equation}\label{ineq div:1}
		u^{h_k}_i(t+h_k)\le \frac12\left( \po(i-j)+\po(i-\ell) +(d_-^{\e,\po}(u^{h_k}(t))_j+u^{h_k}_\ell(t)  \right)+ \frac{C h_k}{R-c_\phi \e},
	\end{equation}  
	for all $i\in \kz\setminus E^{h_k}(t)$. Infimizing in   $j,\ell$ over $E^{h_k}(t)$  in \eqref{ineq div:1} and using again \eqref{def d,sd} and \eqref{sd+=},  we conclude
	\begin{equation}\label{eq:iterk}
		\begin{split}
			u^{h_k}_i(t+h_k)\le  d^{h_k}_i(t) +h_k\frac{C}{R-c_\phi\e_k} \le d^{h_k}_i(t) +h_k\frac{C}{R}.
		\end{split}
	\end{equation}
	provided $h_k,\e_k$ are small enough depending on $R$, and for a possibly larger value of $C$.
	As a consequence of~\eqref{eq:iterk}, we obtain
	\begin{equation}\label{numero}
		-D^*_{\e_k} z^{h_k}(t+h_k)\le   \frac{C}{R}  \qquad\text{in }\{ \widehat d^{h_k}(\cdot,t)\geq R\}\cap \kz.
	\end{equation}
	Using again Lemma \ref{lemma d star to div} and the convergences of $E_{h_k}$ and $d_{h_k}$ it follows that 
	\[ \Div z \le \frac{C}{R} \qquad\text{in }\{(x,t)\in\R^N\times (0, T^*):\, d(x,t)>R\}
	\]
	in the sense of distributions. Hence $\Div z$  is a Radon measure in $\R^N\times (0,T^*)\setminus E$, and 
	$(\Div z)^+\in L^\infty(\{(x,t)\in\R^N\times (0,T^*):\, d(x,t)\geq\delta\})$ for every $\delta>0$.
	
	On the other hand, note that for every  $i\in \kz$  it holds
	\[
	d^{h_k}(t)\geq d^{h_k}_i(t)-\po(\cdot-i).
	\]
	Thus, by Lemma \ref{speriamo che sia vero} and by comparison as before,  we get
	\[
	u^{h_k}_i(t+h_k)\geq d^{h_k}_i(t)-\phi^{h_k}_0 = d^{h_k}_i(t)-(C+1)\sqrt{h_k}.
	\]
	Combining the above inequality with \eqref{eq:iterk},  we deduce for all $t\in (0, T^*)\setminus\mathcal{N}$ and any $\delta>0$ that
	\[
	\sup_{\{\widehat d_{h_k}(\cdot,t)\geq\delta\}\cap \kz}|u^{h_k}(t+h_k)-d^{h_k}(t)|\leq \sqrt{h_k}(C+2),
	\]
	provided that $k$ is large enough.  In particular, recalling also \eqref{conv1},  we deduce that 
	\begin{equation}\label{elleuno}
		\sup_{\kz\cap K}|u^{h_k}(t)-d(\cdot, t)|\to 0\text{ as $k\to\infty$ for $t\in (0, T^*)\setminus \mathcal N$ and for any compact $K\subseteq \R^N\setminus E(t)$},
	\end{equation}
	also with the sequence $\{u^{h_k}\}$ locally (in space and time) uniformly bounded.

	\noindent\textbf{Step 3:} (The subdifferential inclusion). 
	It remains to show that
    \begin{equation}\label{conv subdiff}
		z\in \pa \phi(\nabla d)\quad \text{a.e. in }\R^N\times (0,T^*)\setminus E.
	\end{equation}
	Recall that $ \xi \in \pa \phi(\eta)$ if and only if $ \xi\in \{ v: \po(v)\le 1,\ v\cdot \eta \ge \phi(\eta) \}$.  Since one inequality has been proved in \eqref{subdiff 0}, we show the other one. Consider a test function $\eta\ge 0,\ \eta\in C^{\infty}_c((\R^N\times(0,T^*))\setminus E)$. Let $\sigma>0$ and set $d_\sigma\in C^\infty(\R^N\times (0,T^*))$ as $d_\sigma=d\ast \rho_\sigma$, where $\rho_\sigma$ are space-time mollifiers. 
	Obviously 
	\begin{equation}\label{eq rewr}
		\begin{split}
			\sum_{i,j\in \kz }  z^{h_k}_{ij}(t)\eta(i,t) (u^{h_k}_i(t)&-u^{h_k}_j(t))=  \sum_{i,j\in \kz }z^{h_k}_{ij}(t)\eta(i,t)) (d_\sigma(i,t)-d_\sigma(j,t)) \\
			&+ \sum_{i,j\in \kz } z^{h_k}_{ij}(t) \eta(i,t) \left(u_i^{h_k}(t)-d_\sigma(i,t)-u_j^{h_k}(t)+d_\sigma(j,t) \right).
		\end{split}
	\end{equation}
	In turn,  Lemma \ref{lemma d star to div} implies that
	\begin{equation}\label{eq rewr 1}
		\lim_{k\to \infty} \e_k^{N}\int \left( \sum_{i,j\in \kz }z^{h_k}_{ij}(t)\eta(i,t) \frac{d_\sigma(i,t)-d_\sigma(j,t)}{\e_k}\right)\ud t= \iint  z\cdot \nabla d_\sigma \, \eta\ud x\ud t. 
	\end{equation}
	Let us thus show that
	\begin{align}
		\lim_{\sigma\to 0}  \lim_{k\to \infty}  \e_k^{N}  \int \sum_{i,j\in \kz }  \left( z^{h_k}_{ij}(t) \eta(i,t)  \frac{u^{h_k}_i(t)-d_\sigma(i,t)-u^{h_k}_j(t)-d_\sigma(j,t)}{\e_k} \right)\ud t=0, \label{error term1}
	\end{align} 
	We set for every $t\in (0,T_h^*)$ and $\sigma>0$
	\begin{align*}
		&m_{k,\sigma}(t):= \min_{i\in \textup{supp}(\eta)\cap \kz} (u^{h_k}_i(t)-d_\sigma(i,t)),\\
		&M_{k,\sigma}(t):= \max_{i\in \textup{supp}(\eta)\cap \kz} (u^{h_k}_i(t)-d_\sigma(i,t)).
	\end{align*}
	The convergence \eqref{elleuno} implies that these quantities are uniformly bounded and 
	\begin{equation}\label{mek}
		\lim_{\sigma\to 0}\lim_{k\to +\infty} m_{k,\sigma}(t)=0,\qquad  \lim_{\sigma\to 0}\lim_{k\to +\infty} M_{k,\sigma}(t)=0, 
	\end{equation}
	uniformly for all $t\notin \mathcal N$. For all times $t\in (0,T^*)\setminus \mathcal N$  it holds 
	\begin{align}
		&\e_k^N\sum_{i,j\in \kz } z^{h_k}_{ij}(t) \eta(i,t) \frac{u_i^{h_k}(t)-d_\sigma(i,t)-u^{h_k}_j(t)+d_\sigma(j,t)}{\e_k} \nonumber\\
		&=\e_k^N\sum_{i,j \in \kz} z^{h_k}_{ij}(t) \eta(i,t) \frac{(u_i^{h_k}(t)-d_\sigma(i,t)-m_{k,\sigma}(t))-(u^{h_k}_j(t)-d_\sigma(j,t)-m_{k,\sigma}(t))}{\e_k}\nonumber\\
		&=\e_k^N\sum_{i\in \kz} (u^{h_k}_i(t)-d_\sigma(i,t)-m_{k,\sigma}(t))\sum_{j\in \kz}\left(  \frac{z_{ij}^{h_k}(t)-z_{ji}^{h_k}(t)}{\e_k}\eta(i,t)+ z_{ji}^{h_k}(t) \frac{\eta(i,t)-\eta(j,t)}{\e_k} \right)\label{eq err 1}.
	\end{align} 
	For $k$ large enough, since the support of $\eta$ is at positive distance from $E$, by the bound \eqref{numero} one has $D^*_{\e_k} z^{h_k}(t)\ge - c(\delta)$  on the support for $h_k$ small enough. Thus  it holds 
	\begin{equation*}
		\begin{split}
			\e_k^N \sum_{i\in \kz} (u^{h_k}_i(t)-d_\sigma(i,t)-m_{k,\sigma}(t))\eta(i,t) \sum_{j\in\kz }\frac{z_{ij}^{h_k}(t)-z_{ji}^{h_k}(t)}{\e_k}\\
			\ge   -c(\delta) \e_k^N\sum_{i\in\kz} (u_i^{h_k}(t)-d_\sigma(i,t)-m_{k,\sigma}(t))\eta(i,t).
		\end{split}
	\end{equation*}
	Recalling that $\#(\textup{supp}(\eta)\cap \kz)=O(h_k^{-N})$ uniformly in time, by uniform convergence and \eqref{elleuno} we conclude that 
	\begin{equation}\label{eq err 2}
		\lim_{\sigma \to 0}\liminf_{k\to \infty}\e_k^N \int \sum_{i\in \kz} (u^{h_k}_i(t)-d_\sigma(i,t)-m_{\e,k}(t))\eta(i,t) \sum_{j\in\kz }\frac{z^{h_k}_{ij}(t) -z_{ji}^{h_k}(t)}{\e_k}\ud t\ge 0.
	\end{equation}
	The other term in \eqref{eq err 1} can be estimated using the Lipschitz constant of $\eta$: 
	\begin{equation*}
		\begin{split}
			\left\lvert\int  \sum_{i,j\in \kz}  \e_k^N(u_i^{h_k}(t)-d_\sigma(i,t)-m_{\e,k}(t)) z_{ji}^{h_k}(t)\frac{\eta(i,t)-\eta(j,t)}{\e_k} \ud t\right\rvert\\
			\le \|\nabla\eta\|_\infty \e_k^N\int \sum_{i,j\in \kz} (u_i^{h_k}(t)-d_\sigma(i,t)-m_{\e,k}(t)) \alpha^{h_k}_{ji}\frac{|i-j|}{\e_k}\ud t\to 0
		\end{split}
	\end{equation*}
	letting first $k\to +\infty$ and then $\sigma\to 0$, thanks to \eqref{elleuno} and \eqref{mek}. Note now that adding and subtracting $M_{\e,k}(t)$ to \eqref{error term1} instead of $m_{\e,k}(t)$ and reasoning as above, one proves  that
	\begin{equation}\label{eq err 4}
		\begin{split} 
			&\lim_{\sigma\to 0}\limsup_{k\to \infty}\e_k^N\int \left(  \sum_{i\in \kz} ((u^{h_k}_i(t)-d_\sigma(i,t)-M_{\e,k}(t))\eta(i,t) \sum_{j\in\kz }\frac{z^{h_k}_{ij}(t) -z^{h_k}_{ji}(t)}{\e_k}\right)\ud t\le 0, \\
			&\lim_{\e\to 0}\lim_{k\to \infty}\e_k^N  \int \left\lvert \sum_{i,j\in \kz}  ((u^{h_k}_i(t)-d_\sigma(i,t)-M_{\e,k}(t)) z_{ji}^{h_k}(t)\frac{\eta(i,t)-\eta(j,t)}{\e_k}\right\rvert \ud t = 0.
		\end{split}
	\end{equation}
	Combining \eqref{eq err 1}, \eqref{eq err 2} and \eqref{eq err 4}, we conclude   \eqref{error term1}. 
	
	Integrating in time \eqref{eq rewr}  and combining \eqref{eq rewr 1} and \eqref{error term1}, since $\nabla d_\sigma=\rho_\sigma \ast\nabla d\to \nabla d$ pointwise a.e. and are uniformly bounded in $L^\infty(\R^N\times (0,T^*);\R^N)$, it holds 
	\begin{equation*}\label{eq conv 1}
		\lim_{k\to \infty}  \e_k^N\int  \left(\sum_{i,j\in \kz }   \eta(i,t)\,z^{h_k}_{ij}(t)  \frac{u^{h_k}_i(t)-u^{h_k}_j(t)}{\e_k} \right)   \ud t=\iint z\cdot \nabla d \,\eta \ud x\ud t.
	\end{equation*}
	The convergence above can be paired with the lower semicontinuity  of the $\Gamma$-convergence of the discrete total variations (which follows from an adaptation of classical arguments, see e.g. \cite{ChaKre}) and $z\h_{ij}(u\h_i-u\h_j)=\alpha\h_{ij}|u\h_i-u\h_j|$ to obtain
	\begin{align*}
		\iint \phi(\nabla d)\eta &\le \liminf_{k\to \infty} \e_k^N \int \left(\sum_{i,j\in \kz } \eta(i,t)\,\alpha^{h_k}_{ij} \frac{|u^{h_k}_i(t)-u^{h_k}_j(t)|}{\e_k} \right)\ud t\\
		&=\liminf_{k\to \infty} \e_k^N\int  \left(\sum_{i,j\in \kz }   \eta(i,t)\,z^{h_k}_{ij}(t)  \frac{u^{h_k}_i(t)-u^{h_k}_j(t)}{\e_k} \right)  \ud t =\iint z\cdot \nabla d\, \eta, 
	\end{align*}
	which shows that $\phi(\nabla d)=z\cdot \nabla d$ a.e. on the support of $\eta$, from which we deduce \eqref{conv subdiff}.
\end{proof}
We conclude this section by observing that the discrete scheme converges to the unique weak flow (in the sense of Definition~\ref{Defsol}) starting from $E_0$ for ``generic'' initial data $E_0$, i.e. whenever fattening does not occur. More precisely, we have the following Corollary.
\begin{corollary}Let $u_0\in \textup{UC}(\R^N)$ and for every $\lambda\in \R$ let $\overline E^h_\lambda$ be the closed space-time tube of the $h$-discrete evolution starting from $\{u_0\le\lambda\}$; i.e., as in \eqref{tubidiscreti} with $E_0=\{u_0\le\lambda\}$. Then, there exists a countable set $\mathcal N$ such that for all $\lambda\in \R^N\setminus\mathcal N$
	\[
	\overline E^{h}_\lambda\stackrel{\mathcal K}{\longrightarrow} E_\lambda\quad \text{ in }\R^N\times[0,+\infty)
	\]
	as $h\to 0$, where $E_\lambda$ is the unique weak flow in the sense of Definition~\ref{Defsol} starting from $\{u_0\le\lambda\}$.
\end{corollary}
\begin{proof}
	It follows by combining Theorems~\ref{themthm}  and \ref{th:phiregularlevelset}.
\end{proof}

\section{Numerical experiments}\label{sec:Num}

\begin{figure}[htb]
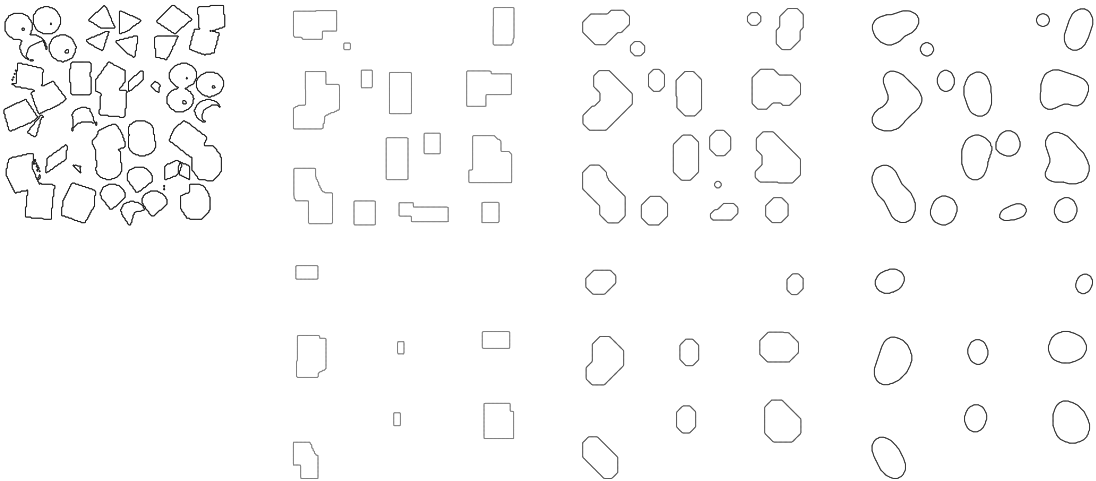

	\resizebox{0.2\textwidth}{!}{\input{figures/sh_init}}\hspace{.05\textwidth} 
	\resizebox{.2\textwidth}{!}{\input{figures/sh_squ50}}\hspace{.05\textwidth} 
	\resizebox{.2\textwidth}{!}{\input{figures/sh_oct50}}\hspace{.05\textwidth} 
	\resizebox{.2\textwidth}{!}{\input{figures/sh_iso50}}\\[5mm]
	\hspace{.2\textwidth}\hspace{.05\textwidth} 
	\resizebox{.2\textwidth}{!}{\input{figures/sh_squ120}}\hspace{.05\textwidth} 
	\resizebox{.2\textwidth}{!}{\input{figures/sh_oct120}}\hspace{.05\textwidth} 
	\resizebox{.2\textwidth}{!}{\input{figures/sh_iso120}}
	\caption{An initial datum and evolutions for square, octagonal and ``almost isotropic''
		anisotropies, at two different times.}\label{fig:shapes}
\end{figure}
\begin{figure}[htb]
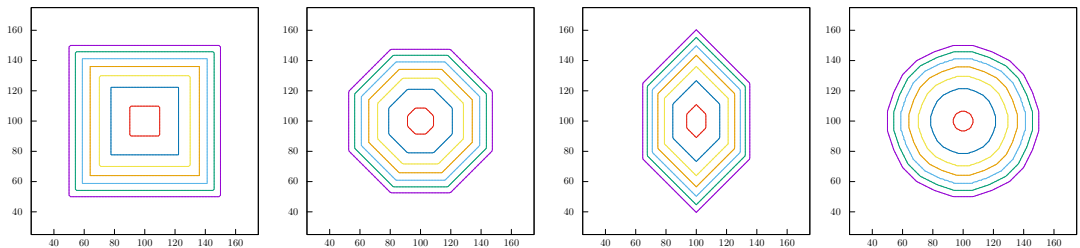

	\resizebox{0.23\textwidth}{!}{\input{figures/squ50s}}
	\ \resizebox{0.23\textwidth}{!}{\input{figures/oct50s}}
	\ \resizebox{0.23\textwidth}{!}{\input{figures/dia50s}}\ \resizebox{0.23\textwidth}{!}{\input{figures/iso50s}}
	\caption{Wulff shapes of initial radius $R_0=50$ evolved at times $t=0$, $200$, $400,\dots,1200$ for four different 
		anisotropies (square, octagon, diamond and ``almost isotropic'').}\label{fig:Wulffshapes}
\end{figure}

We  show some numerical experiments to illustrate our results, in dimension~2 (an implementation
in 3D is currently being developed). We follow the implementation
described in~\cite{ChaDar} (see also~\cite{DarbonChambolle12}), except that now the distance
is properly computed using using the inf/sup-convolution formulas~\eqref{def d,sd}.
The (exact) numerical resolution of the discrete ROF functional is computed using
Hochbaum's parametric maximum flow algorithm~\cite{Hochbaum01,Hochbaum13}, 
implemented upon the maxflow/mincut implementation
of Boykov and Kolmogorov~\cite{BoykovKolmogorov-maxflow}. 
{This particular algorithm has the advantage
to provide an \textit{exact} solution of the ROF problem, up to computer precision.}
Other implementations of the
algorithm
yielding approximate minimizers have been considered for instance in~\cite{Chambolle04,Obermanetal2011},
of course they work in practice and allow to address more (an)isotropies than the current
work, yet the joint convergence as $\e=h\to 0$ is not clear in these contexts.
For numerical speedup, the infimum and supremum of definition \eqref{def d,sd} are computed only in a neighborhood of fixed size and not on the whole grid. 
{We expect this to yield, in general, an error of order $C\e$ with $C$ getting smaller as the width of the strip increases, however 
we} observe that Corollary~\ref{cor neigh} in Appendix~\ref{app:finite} justifies this restriction {(showing that $C=0$)} in some cases, notably the case $\phi=\|\cdot\|_{\ell^1}, \po=\|\cdot\|_{\ell^\infty}$, {see Fig.~\ref{fig:Wulffshapes}, left, for which the sup/inf are
in fact min/max which are reached very close to the evolving boundary (as one can chose $\ell_1=1$  in Lemma \ref{lemma neigh}).}
Similarly, the ROF minimization is only performed in a neighborhood of the boundary {(one can show that this does not affect the solution in a smaller neighborhood, hence the overall error is the same as when computing the distance in a strip only)}. 

The code is available at \url{https://plmlab.math.cnrs.fr/chambolle/discretecrystals/} (implemented
in \texttt{C/C++} and running on GNU/linux with \texttt{gcc}).

Figure~\ref{fig:shapes} shows three examples of flows from the same starting set, composed of 
random shapes. The anisotropies are square (nearest neighbours interactions), octagonal (next nearest neighbours, weighted so that the corresponding Wulff shape is a regular octagon), and ``almost
isotropic'', which is generated by the interactions in the $12$ directions $e\in \{(0,\pm 1), (\pm 1,0), (\pm 1,\pm 2), (\pm 1,\pm 3)\}$ weighted so that the Wulff shape is a polygon with 24 facets of equal lengths. 
{This is obtained by setting the weights in
the discrete total variation as $.131/\text{length}(e)$
for each direction $e$, so that the total perimeter
of the unit Wulff shape is $24\times (2\times.131)\approx 2\pi$, in the hope that the
corresponding crystalline curvature will be close
to the Euclidean one.}
\begin{figure}[htb]
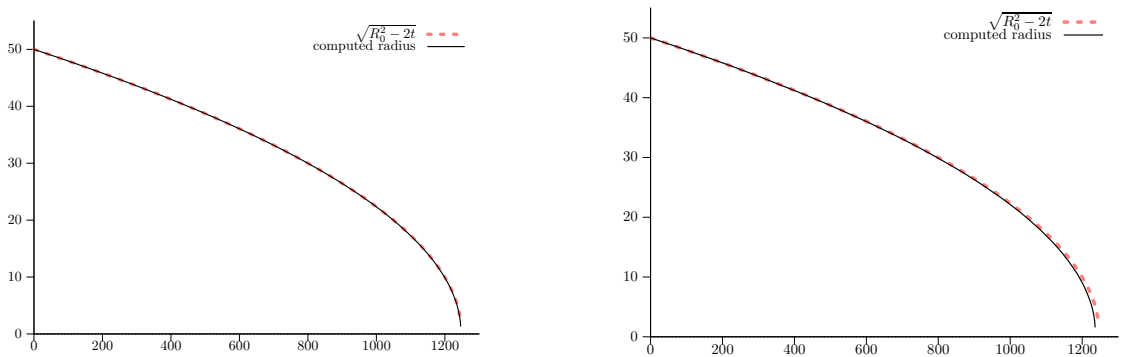

	\resizebox{0.45\textwidth}{!}{\input{figures/radii_squ_50}}\hfill
	\resizebox{0.45\textwidth}{!}{\input{figures/radii_oct_50}}
	\caption{Evolution of the radius for the square (left)
		and octogonal (right) anisotropies.}\label{fig:radiiSqOc}
\end{figure}
\begin{figure}[htb]
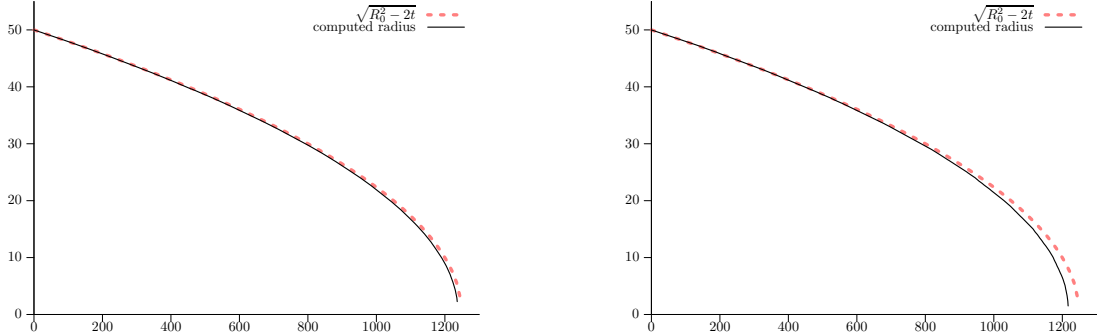

	\resizebox{0.45\textwidth}{!}{\input{figures/radii_dia_50}}\hfill
	\resizebox{0.45\textwidth}{!}{\input{figures/radii_iso_50}}
	\caption{Evolution of the radius for the diamond (left)
		and ``almost isotropic'' (right) anisotropies.}\label{fig:radiiDiIs}
\end{figure}

Then, we estimate the decay of the radius of an initial Wulff shape $\W_{R_0}=\{\phi \le R_0\}$ along the evolution,
up to extinction. In our experiment, $R_0=50$. It is well known that the solution
is the Wulff shape of radius $R(t)=\sqrt{R_0^2-2(N-1)t}$ (where here $N=2$). The evolutions
are depicted in Figure~\ref{fig:Wulffshapes}. We use the same anisotropies as in figure~\ref{fig:shapes},
with additionally a ``diamond'' Wulff shape generated by the directions $(0,\pm 1)$, $(\pm 1,\pm 2)$
and with sides of equal lengths. In all cases, the weights have been calibrated so that the perimeters
of the Wulff shapes are $6.28\approx 2\pi$.

The plots in Figure~\ref{fig:radiiDiIs} show that the decay of the radii is remarkably close
to the theoretical prediction, even if this is less precise when more directions of interactions are involved, near extinction. This might be due in part to the fact that the computation of the distance
through truncated variants of~\eqref{def d,sd} become less precise.

\begin{figure}[htb]
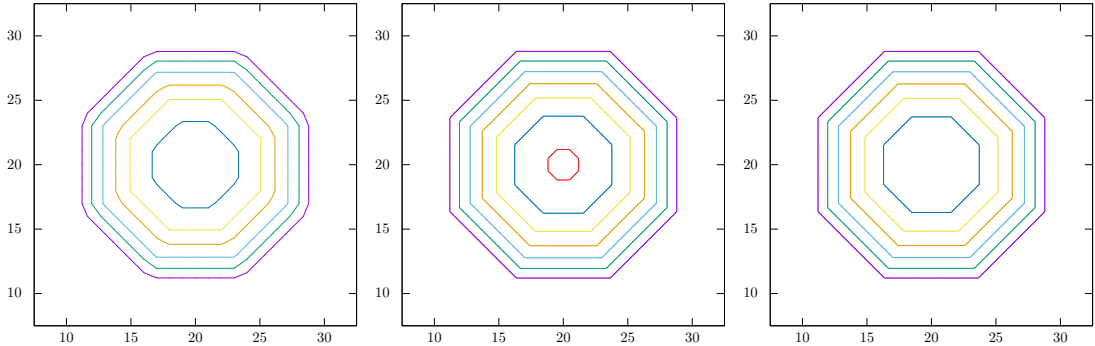

	\resizebox{0.32\textwidth}{!}{\input{figures/octos40}}\ 
	\resizebox{0.32\textwidth}{!}{\input{figures/octos400}}\ 
	\resizebox{0.32\textwidth}{!}{\input{figures/octos400_50}}
	\caption{Evolution of an initial octagon 
		with $R_0=10$ at times $0,7,14,\dots$. Left: $\e=1$, $h=0.1$, middle: $\e=0.1$, $h=0.1$, right: $\e=0.1$, $h=0.5$.}\label{fig:octos}
\end{figure}
\begin{figure}[htb]
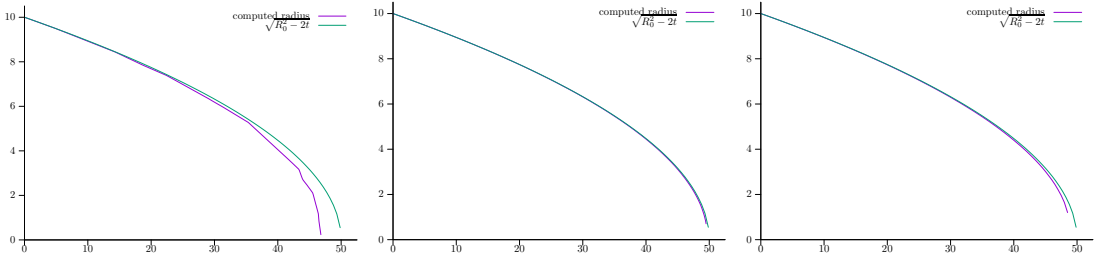

	\resizebox{0.32\textwidth}{!}{\input{figures/rad_oct_40}}\ 
	\resizebox{0.32\textwidth}{!}{\input{figures/rad_oct_400}}\ 
	\resizebox{0.32\textwidth}{!}{\begin{tikzpicture}[gnuplot]
\gpcolor{color=gp lt color border}
\gpsetlinetype{gp lt border}
\gpsetdashtype{gp dt solid}
\gpsetlinewidth{1.00}
\draw[gp path] (1.012,0.616)--(0.832,0.616);
\node[gp node right] at (0.828,0.616) {$0$};
\draw[gp path] (1.012,2.106)--(0.832,2.106);
\node[gp node right] at (0.828,2.106) {$2$};
\draw[gp path] (1.012,3.597)--(0.832,3.597);
\node[gp node right] at (0.828,3.597) {$4$};
\draw[gp path] (1.012,5.087)--(0.832,5.087);
\node[gp node right] at (0.828,5.087) {$6$};
\draw[gp path] (1.012,6.578)--(0.832,6.578);
\node[gp node right] at (0.828,6.578) {$8$};
\draw[gp path] (1.012,8.068)--(0.832,8.068);
\node[gp node right] at (0.828,8.068) {$10$};
\draw[gp path] (1.012,0.616)--(1.012,0.436);
\node[gp node center] at (1.012,0.308) {$0$};
\draw[gp path] (3.095,0.616)--(3.095,0.436);
\node[gp node center] at (3.095,0.308) {$10$};
\draw[gp path] (5.178,0.616)--(5.178,0.436);
\node[gp node center] at (5.178,0.308) {$20$};
\draw[gp path] (7.261,0.616)--(7.261,0.436);
\node[gp node center] at (7.261,0.308) {$30$};
\draw[gp path] (9.343,0.616)--(9.343,0.436);
\node[gp node center] at (9.343,0.308) {$40$};
\draw[gp path] (11.426,0.616)--(11.426,0.436);
\node[gp node center] at (11.426,0.308) {$50$};
\draw[gp path] (1.012,8.441)--(1.012,0.616)--(11.947,0.616);
\node[gp node right] at (10.479,8.107) {computed radius}; 
\gpcolor{rgb color={0.580,0.000,0.827}}
\draw[gp path] (10.663,8.107)--(11.579,8.107);
\draw[gp path] (1.012,8.068)--(1.116,8.031)--(1.220,7.993)--(1.324,7.955)--(1.429,7.917)%
  --(1.533,7.879)--(1.637,7.840)--(1.741,7.802)--(1.845,7.763)--(1.949,7.723)--(2.053,7.684)%
  --(2.158,7.644)--(2.262,7.604)--(2.366,7.564)--(2.470,7.524)--(2.574,7.484)--(2.678,7.443)%
  --(2.782,7.402)--(2.887,7.361)--(2.991,7.320)--(3.095,7.278)--(3.199,7.236)--(3.303,7.194)%
  --(3.407,7.151)--(3.511,7.108)--(3.616,7.065)--(3.720,7.021)--(3.824,6.978)--(3.928,6.934)%
  --(4.032,6.890)--(4.136,6.845)--(4.240,6.800)--(4.345,6.755)--(4.449,6.709)--(4.553,6.663)%
  --(4.657,6.616)--(4.761,6.570)--(4.865,6.523)--(4.969,6.476)--(5.074,6.428)--(5.178,6.380)%
  --(5.282,6.331)--(5.386,6.282)--(5.490,6.232)--(5.594,6.182)--(5.698,6.131)--(5.803,6.080)%
  --(5.907,6.029)--(6.011,5.977)--(6.115,5.925)--(6.219,5.872)--(6.323,5.818)--(6.427,5.764)%
  --(6.532,5.709)--(6.636,5.654)--(6.740,5.598)--(6.844,5.542)--(6.948,5.485)--(7.052,5.427)%
  --(7.156,5.369)--(7.261,5.309)--(7.365,5.248)--(7.469,5.187)--(7.573,5.126)--(7.677,5.063)%
  --(7.781,5.000)--(7.885,4.935)--(7.990,4.870)--(8.094,4.803)--(8.198,4.735)--(8.302,4.667)%
  --(8.406,4.597)--(8.510,4.526)--(8.614,4.452)--(8.719,4.377)--(8.823,4.301)--(8.927,4.223)%
  --(9.031,4.144)--(9.135,4.062)--(9.239,3.979)--(9.343,3.893)--(9.448,3.805)--(9.552,3.716)%
  --(9.656,3.624)--(9.760,3.530)--(9.864,3.429)--(9.968,3.326)--(10.072,3.220)--(10.177,3.107)%
  --(10.281,2.990)--(10.385,2.864)--(10.489,2.729)--(10.593,2.590)--(10.697,2.437)--(10.801,2.262)%
  --(10.906,2.070)--(11.010,1.839)--(11.114,1.510);
\gpcolor{color=gp lt color border}
\node[gp node right] at (10.479,7.799) {$\sqrt{R_0^2-2t}$}; 
\gpcolor{rgb color={0.000,0.620,0.451}}
\draw[gp path] (10.663,7.799)--(11.579,7.799);
\draw[gp path] (1.012,8.068)--(1.122,8.029)--(1.233,7.989)--(1.343,7.949)--(1.454,7.909)%
  --(1.564,7.868)--(1.675,7.827)--(1.785,7.786)--(1.896,7.745)--(2.006,7.704)--(2.117,7.662)%
  --(2.227,7.620)--(2.337,7.578)--(2.448,7.536)--(2.558,7.493)--(2.669,7.450)--(2.779,7.407)%
  --(2.890,7.363)--(3.000,7.319)--(3.111,7.275)--(3.221,7.231)--(3.332,7.186)--(3.442,7.141)%
  --(3.552,7.096)--(3.663,7.050)--(3.773,7.004)--(3.884,6.958)--(3.994,6.912)--(4.105,6.865)%
  --(4.215,6.817)--(4.326,6.770)--(4.436,6.722)--(4.547,6.673)--(4.657,6.624)--(4.767,6.575)%
  --(4.878,6.525)--(4.988,6.475)--(5.099,6.425)--(5.209,6.374)--(5.320,6.323)--(5.430,6.271)%
  --(5.541,6.218)--(5.651,6.166)--(5.762,6.112)--(5.872,6.058)--(5.982,6.004)--(6.093,5.949)%
  --(6.203,5.894)--(6.314,5.838)--(6.424,5.781)--(6.535,5.723)--(6.645,5.665)--(6.756,5.607)%
  --(6.866,5.547)--(6.977,5.487)--(7.087,5.426)--(7.197,5.365)--(7.308,5.302)--(7.418,5.239)%
  --(7.529,5.175)--(7.639,5.110)--(7.750,5.044)--(7.860,4.977)--(7.971,4.909)--(8.081,4.840)%
  --(8.192,4.769)--(8.302,4.698)--(8.412,4.625)--(8.523,4.551)--(8.633,4.475)--(8.744,4.398)%
  --(8.854,4.320)--(8.965,4.239)--(9.075,4.157)--(9.186,4.073)--(9.296,3.986)--(9.407,3.898)%
  --(9.517,3.807)--(9.627,3.713)--(9.738,3.617)--(9.848,3.517)--(9.959,3.413)--(10.069,3.306)%
  --(10.180,3.194)--(10.290,3.077)--(10.401,2.955)--(10.511,2.825)--(10.622,2.688)--(10.732,2.540)%
  --(10.842,2.381)--(10.953,2.205)--(11.063,2.007)--(11.174,1.776)--(11.284,1.486)--(11.395,1.026);
\gpcolor{color=gp lt color border}
\draw[gp path] (1.012,8.441)--(1.012,0.616)--(11.947,0.616);
\gpdefrectangularnode{gp plot 1}{\pgfpoint{1.012cm}{0.616cm}}{\pgfpoint{11.947cm}{8.441cm}}
\end{tikzpicture}
}
	\caption{Evolution of the radius for an initial octagon 
		with $R_0=10$ until the vanishing time $t=50$. Left: $\e=1$, $h=0.1$, middle: $\e=0.1$, $h=0.1$, right: $\e=0.1$, $h=0.5$.}\label{fig:radiiOcto}
\end{figure}

Finally, we perform the same experiment with varying $\e$ and $h$.
We observe that the results look remarkably close even if, at low resolution,
the error becomes huge when the size of the Wulff shape is of the order of
the discretization. Figure~\ref{fig:octos} shows the shapes. Observe that
the shape at time $t=49$ is only computed for $\e=0.1$ and $h=0.1$
(the shape vanishes before for the two other experiments). On the other hand,
this computation took more than one hour, while the case $\e=1$ took less than
a minute and the case $\e=0.1$, $h=0.5$ a bit less than an hour.
Figure~\ref{fig:radiiOcto} shows the decay of the radii,
which should be $\sqrt{R_0^2-2t}$ for $R_0=10$ and $t\in [0,50]$.

\appendix

\section{Proof of Lemma \ref{speriamo che sia vero}}\label{sec:proofscsv}

We build here a supersolution to Problem~\eqref{discrete EL resc} when $g=\po$.
Let us first recall some notation and results concerning zonotopes (see e.g. \cite{McM}). 
Recall that $\mathcal E=\{\pm e_k\}_{k=1}^m\subseteq \Z^N$ where, without loss of generality, the vectors $e_1,\dots,e_m$ span the whole $\R^N$. Given a non-negative interaction function $\beta\in X$, we assume that $\beta=0$ on $\Z^N\setminus \mathcal E$ and that $\beta(-i)=\beta(i)$ for every $i\in\Z^N.$
The anisotropy $\phi$ associated to $\beta,$ as defined in \eqref{def phi}, is   such that its $1$-Wulff shape $\W_1\subseteq \R^N$ is a zonotope, which can be expressed as the Minkowski sum 
\[ \W_1=\sum_{e\in\mathcal E} \beta(e) {[-e,e]}=\sum_{k=1}^m 2\beta(e_k){[-e_k,e_k]}, \]
{where ${[-e,e]}\subseteq\R$ denotes the  closed segment from $-e$ to $e$.}
Alternatively, one can define the zonotope $\W_1$ as the image of a cube under an affine map. Indeed, it holds
\begin{equation}\label{def proj cube}
	\W_1=V(Q^{(m)})
\end{equation}
where $V=(2\beta(e_1) e_1, \dots, 2\beta(e_m) e_m)\in\R^{N\times m}$ and $Q^{(m)}={[-1,1]^m}$. Since the set $\mathcal E$ is uniquely defined up to sign changes, the matrix $V$ is also uniquely detemined up to permutations of columns or sign changes. 

Note that by definition of zonotope any element $x\in {\W_\ell }$ for $\ell>0$ can be written as 
\[ x=\ell\sum_{k=1}^m 2 \beta(e_k) \lambda_k e_k,\]
for suitable coefficients $|\lambda_k|\le 1.$  
We note that (the closure of) a facet $F$ (of non-zero dimension) of the zonotope $\W_\ell$ can be
described in the following form:
\begin{equation}\label{formula facet}
F=\ell\sum_{j=1}^r 2\beta(e_{\sigma(j)}) [-e_{\sigma(j)}, e_{\sigma(j)}] + \ell\sum_{j=r+1}^m 2\beta(e_{\sigma(j)})\e_{\sigma(j)} e_{\sigma(j)},
\end{equation}
where $\sigma$ is a permutation of $\{1,\dots, m\}$, $1\le r\le m$ and $|\e_j|=1$. Moreover (see \cite[page 206]{McM} for details)  the vectors $e_{\sigma(1)},\dots, e_{\sigma(r)}$ uniquely identify 
\[ \{ e\in\mathcal E : e\parallel F \}, \]
and $r$ is uniquely defined as the number of vectors in the family $\mathcal E$ which are parallel to the facet $F$. Analogously, any vertex $v$ of the zonotope $\W_\ell$ is of the form 
\begin{equation}\label{formula vertex}
v=\ell \sum_{j=1}^m  2\beta(e_{\sigma(j)})\e_{\sigma(j)} e_{\sigma(j)},
\end{equation}
where $\e_j\in \{\pm 1\}$ for every $j=1,\dots,m$ and $\sigma$ is a permutation of $\{1,\dots, m\}$. Note however that not every point of this form is a vertex of the zonotope.

\begin{lemma}\label{lemma zonotopes}
There exists $\ell_0>0$ such that for every $\e>0$ and every $\ell\ge \ell_0$, if $i\in\ez$ belongs to $\partial \W_{\e\ell}$, then for each $k\in \{1,\dots,m\}$ either one of the following holds:
\begin{itemize}
	\item [i)] neither $i+ \e e_k$ nor $i-\e e_k$ belong to $ \partial \W_{\e \ell}$. In this case it holds either $\po(i+\e e_k)>\po(i)>\po(i-\e e_k)$   or $\po(i-\e e_k)>\po(i)>\po(i+\e e_k)$;
	\item[ii)] one between $i\pm \e e_k$ belongs to $\partial \W_{\e \ell}$. In this case $\po(i\pm \e e_k)\ge \ell$ and  it holds
	\begin{equation}\label{length}
		\#(( i+\e\Z e_k)\cap \partial \W_{\e\ell}   )\ge 2[\ell/\ell_0].
	\end{equation}
\end{itemize}
\end{lemma}

\begin{proof}
By scaling, it suffices to prove the result in the case $\e=1.$ We take $\ell_0$ such that
\begin{equation}\label{choice ell0}
	\ell_0\ge \max_{k=1,\dots, m}\frac 1{2\beta(e_k)}
\end{equation}
and remark that $\ell_0\in (0,+\infty).$ Note that the choice \eqref{choice ell0} implies for every $j=1,\dots, m$ that
\[ 
|[-2\ell \beta(e_j) e_j ,2\ell \beta(e_j) e_j]|=4 \ell \beta(e_j) |e_j|\ge 2   \frac\ell{\ell_0}|e_j|. 
\]
We then fix $i\in\partial \W_\ell\cap \Z^N$ and $e_k\in\mathcal E$. We have to distinguish two cases.\\
\textbf{Case 1.} There exists a facet $F\ni i$ of $\W_\ell$ such that $e_k\parallel F$. By \eqref{formula facet} we then see that
\[ i\in 2\ell \beta(e_k)[-e_k,e_k]+j , \]
where $j\in F$. This implies in particular that  $\{ n\in\Z : i+n e_k\in F \} $     is an interval of $\Z$ containing $0$. Furthermore, by the assumption \eqref{choice ell0}, it contains at least $[2\ell|e_k|/\ell_0]$ points and we conclude \eqref{length}. Since $i$ and one between $i\pm e_k$ belong to $\partial \W_\ell$, then $\po(i\pm e_k)\ge \ell$ by convexity. \\
\textbf{Case 2.} For every facet $F\ni i$ of $\W_\ell$ it holds $e_k\nparallel F$. Let us fix a facet $F \ni i$ and note that
by \eqref{formula facet} and  up to relabelling the indexes, it holds
\[  
i\in \ell \sum_{j=1}^r 2\beta(e_j) [-e_j,e_j] +  \ell \sum_{j=r+1}^m 2\beta(e_j)\e_j e_j , 
\]
with $k>r$ and $|\e_j|=1$ for $j=r+1,\dots, m$. Recalling \eqref{def proj cube}, we see that 
\[i-\e_k e_k = \ell V (y - \frac{\e_k}{\ell \beta(e_k)}\tilde e_k),\]
where  $\tilde e_1, \dots \tilde e_m$ denotes the canonical base of $\R^m$ and  $y\in \sum_{j=1}^r  [-\tilde e_j,\tilde e_j] +  \sum_{j=r+1}^m \e_j \tilde e_j\subseteq \partial Q^{(m)}$. By the choice \eqref{choice ell0} and since $k>r$, one deduces that $y - \frac{\e_k}{\ell \beta(e_k)}\tilde e_k\in Q^{(m)}$, thus $i-\e_k e_k\in \overline \W_\ell$. Since then $e_k\nparallel F$ for any facet containing $i$, it must hold $\po(i-\e_k e_k)<\ell$. By convexity one easily concludes that $\po(i+\e_k e_k)>\ell$, which shows i).
\end{proof}

We now define a calibration $z_{ij}$ for every $(i,j)\in \left(\{ \po > \e \ell_0 \}\cap \ez\right)\times \ez$. Fix  $i\in \ez$ with $\po(i)>\e \ell_0$.  In the following we write $i\sim j$ if $\frac{i-j}\e\in\mathcal E$. We start defining
\begin{equation}
z_{ij}=\begin{cases}\label{def z 1}
	0 &\text{if } j\not\sim i\\
	-\beta(e_k) &\text{if }j=i \pm  \e   e_k \text{ and }\po(j)>\po(i)\\
	\beta(e_k) &\text{if }j=i\pm \e  e_k\text{ and }\po(j)<\po(i).
\end{cases}
\end{equation}
In particular, this definition covers case i) in Lemma \ref{lemma zonotopes}. Assume then that there exists $j\sim i$ with $\po(j)=\po(i)$ and $\frac{j-i}\e=e_k\in\mathcal E$. Since $i\in\ez$ and $e_k\in\mathcal E$ fall in case ii) of Lemma \ref{lemma zonotopes}, there exists an interval $[-\underline{n},\bar{n}]\cap \Z$ for $\underline{n},\bar{n}\in \N$ such that
\[( i+\e\Z e_k)\cap \partial W^{\po}_{\po(i)}   = i +  ([-\underline{n},\bar{n}]\cap \Z) \e e_k \]
and moreover
\begin{equation}\label{length2}
\#([-\underline{n},\bar{n}]\cap \Z)\ge 2[\po(i)/(\e\ell_0)].
\end{equation}
Thus, we define $z_{ij}$ as a linear interpolation of the values assumed at the extremal points of $i +  [-\underline{n},\bar{n}] \e e_k $ as
\begin{equation}\label{def z 2}
\begin{split}
	&z_{i+   t \e e_k,i+(t+1)\e e_k}:= \beta(e_k)\left(1-2 \frac{t+\underline{n}+1}{\underline n+ \bar{n}+1}  \right)\quad \forall t\in [-\underline{n}-1,\bar{n}]\cap \Z,\\[1ex]
	&z_{i+   t \e e_k,i+(t-1)\e e_k}:= \beta(e_k)\left(1-2 \frac{-t+\underline{n}+1}{\underline n+ \bar{n}+1}  \right)\quad \forall t\in [-\underline{n},\bar{n}+1]\cap \Z.
\end{split}
\end{equation}
By definition one easily sees that 
\begin{equation}\label{subdiff z}
|z_{ij}|\le \alpha_{ij}^\e,\qquad z_{ij}(\po(i)-\po(j))=\alpha_{ij}^\e|\po(i)-\po(j)|.
\end{equation}
We now show how to bound the divergence $(D^*_\e z)_i$. 
Assume that $\po(i + \e e_k)= \po(i)$ or that $\po(i - \e e_k)= \po(i)$. Then by definition \eqref{def z 2} and by \eqref{length2} one deduces
\begin{equation}
\begin{split}\label{div1}
    &z_{i,i+\e e_k}+z_{i,i-\e e_k} - z_{i+   \e e_k,i} - z_{i- \e e_k,i} =-\frac{4\beta(e_k)}{\underline n+ \bar n+1}\ge -\frac{2\beta(e_k)}{[\po(i)/(\e \ell_0)]} \ge -\frac{C\e }{\po(i)},
\end{split}
\end{equation}
and similarly  if $\po(i - \e e_k)= \po(i)$. If instead $\po(i\pm \e e_k)\neq \po(i)$   and $\po(i\pm \e e_k)\ge \e\ell_0$, one sees that 
\begin{equation}\label{div2}
z_{i,i +\e e_k}+z_{i,i -\e e_k}=0 \ \text{and}\ z_{i +\e e_k,i}+z_{i -\e e_k,i}=0
\end{equation}
Combining \eqref{div1} and \eqref{div2} and recalling \eqref{hp3} we conclude that if $\po(i)\ge \ell_1\e $ then
\begin{equation}\label{bound div big d}
h(D^*_\e z)_i\ge -\frac{c_\phi h}{\po(i)}
\end{equation}
for a suitable positive constant $c_\phi$ depending on $\phi.$

We now illustrate a procedure that allows to extend the calibration above to  $\ez\times \ez$. We set  $C>1$ a sufficiently big constant and define a function $v\in X_\e$ setting
\begin{equation}\label{def v}
v:= \begin{cases}
	\po+\dfrac {C  h}{\po}  \quad &\text{on }\{\po\ge C(\sqrt h\vee \e)\} \cap \ez \\
	 C (\sqrt h\vee \e)+ \frac{h}{\sqrt h\vee \e} \qquad &\text{on } \{\po< C(\sqrt h\vee \e)\}\cap \ez
\end{cases}  .
\end{equation}
A calibration  $w\in Y_\e$ can be defined setting for $i,j\in\ez$
\begin{equation}\label{def z 3}
w_{ij}:=\begin{cases}
	z_{ij}\quad &\text{if } \po(i) \ge 2\sqrt C(\sqrt h\vee \e)\\
	-\alpha_{ij}^\e \quad& \text{if } \po(i) < 2\sqrt C (\sqrt h\vee \e)
\end{cases}.
\end{equation}
Since $x\mapsto x+C h x^{-1}$ is strictly monotone in the region $\{ x\ge \sqrt{Ch} \},$ we can employ  \eqref{subdiff z} to prove that, for every $i,j\in\ez$ with $\po(i)\ge  C (\sqrt h\vee \e)$, it holds
\begin{equation}\label{eq subdiff v}
    w_{ij}(v_i-v_j)=\alpha_{ij}^\e|v_i-v_j|, \quad |w_{ij}|\le \alpha_{ij}^\e.
\end{equation}
Moreover, taking $C$ large enough ensures that whenever $j \sim i$, then 
\begin{equation}
    \begin{split}\label{choice C}
    \po(i)\le 2\sqrt C(\sqrt h\vee \e) &\implies \po(j)\le C(\sqrt h\vee \e)\\
    \po(i)\ge 2\sqrt C(\sqrt h\vee \e) &\implies \po(j)\ge \sqrt C(\sqrt h\vee \e)
\end{split}.
\end{equation}
Thus,   equation \eqref{eq subdiff v} can be directly checked in the case $\po(i)\le 2\sqrt C(\sqrt h\vee \e)$ using the definition \eqref{def z 3}.

Note now that definition \eqref{def z 3} implies $D^*_\e w=0$ in the region $\{\po< 2\sqrt C(\sqrt h\vee \e)\}$ thus we assume  $\po(i)\ge 2\sqrt{C}(\sqrt h\vee \e)$ and estimate $(D^*_\e w)_i$. 
If $\po(i- \e e_k) < 2\sqrt{C} (\sqrt h\vee \e)$ by convexity  $\po(i +\e e_k)>  2\sqrt{C} (\sqrt h\vee \e), $ thus by definition \eqref{def z 3} we get
\[  z_{i,i+\e e_k} - z_{i+\e e_k,i} + z_{i,i-\e e_k} -z_{i-\e e_k,i} = -\beta(e_k)-\beta(e_k) + \beta(e_k) - (-\beta(e_k) )=0. \]
The symmetric case is analogous. On the other hand, if  every $j\sim i$ is in $\{ \po \ge 2\sqrt{C} (\sqrt h\vee \e) \}$ equation \eqref{bound div big d} holds. 
Therefore, we have shown
\begin{equation}\label{bound div w}
hD^*_\e w\ge -\frac{c_\phi h}{\po} \chi_{\{\po \ge \sqrt C(\sqrt h\vee \e)\}}.
\end{equation}
By a direct computation, using \eqref{bound div w} and assuming the $C>c_\phi$,
we see that the pair $(v,w)$ defined above satisfies 
\[ \begin{cases}
hD^*_\e w + v \ge \po\\
w_{ij}(v_i-v_j)=\alpha_{ij}^\e|v_i-v_j|, \quad |w_{ij}|\le \alpha_{ij}^\e.
\end{cases}\]
Recalling the comparison result in Theorem \ref{thm exist}, we  conclude that the solution $u$ to \eqref{discrete EL}  satisfies $u\le v$ in $\ez$.

\section{A remark on the inf/sup-convolution formulas~\eqref{def d,sd}}\label{app:finite}
In this section we show that
in some particular cases, the $\inf,\sup$ in the definition \eqref{def d,sd} can be replaced by $\min,\max$ and that this minimization/maximization procedure can be made in a fixed neighborhood of the point considered.
Yet, our proof also shows that this neighborhood can become very large, depending on the weights of
the interaction, and it seems that we cannot expect in general cases that the $\min,\max$ are actually
reached.

{We introduce the following condition.} There exists $\ell_\phi>0$ such that for every $\e_k\in\{0,\pm 1\}$ for $k=1,\dots, m$, there exists $\ell\le \ell_\phi$ such that 
\begin{equation}\label{rational zonotope}
	\ell \sum_{k=1}^m 2\beta(e_k) \e_k e_k\in \Z^N. 
\end{equation}
Note that this condition is satisfied if {$\beta(e_k)/\beta(e_{k'})\in\Q$} for all $k,k'=1,\dots, m$.
\begin{lemma}\label{lemma neigh}
	There exists $\ell_1>0$ with the following property. For any $i\in \ez$ with $\po(i)\ge\e\ell_1 $ there exists $j\in \ez\setminus \{0\}$ with $\po(j)< \po(i)$ and satisfying
	\begin{equation}\label{additive s d1}
		\po(i)\ge \po(j)+\po(i-j)-c_\phi\e.
	\end{equation}
	If  \eqref{rational zonotope} holds, for any $i\in\ez$ with $\po(i)\ge 2\e\ell_1$ there exists $j\in (\W_{\e\ell_1}\setminus\{0\} ) \cap \ez$ such that 
	\begin{equation}\label{additive sd}
		\po(i)= \po(j)+\po(i-j).
	\end{equation}
	Moreover, for every $R\in (2\e\ell_1,\po(i))$ there exists $j\in \W_{R}\setminus\W_{R-2\e\ell_1}$ such that \eqref{additive sd} holds.
	
\end{lemma}

\begin{proof}
	By scaling we prove the result in the case $\e=1$. 
	Given $i\in\Z^N\setminus \{0\}$, inequality \eqref{additive s d1} follows easily choosing $\ell_1\ge2$, considering $\sigma i\in \R^N\setminus\{0\}$ for an appropriate $\sigma \in (0,1)$ and $j\in\Z^N$ so that $\sigma i\in (j+[0,1]^N)$.

	We now assume \eqref{rational zonotope} and denote by $\ell_\phi$ the radius associated to $\phi$. We then choose $\ell_1=\ell_\phi$. Let us fix $i\in\Z^N$ with $\po(i)=\ell\ge 2 \ell_1$. By \eqref{formula facet} there exist $r>0$, $\e_k, \lambda_k$ with $|\e_k|=1$ and $|\lambda_k|<1$ such that
	\[ i=\ell\left( \sum_{k=1}^r 2\beta(e_k) \e_k  e_k +\sum_{k=r+1}^m \lambda_k 2\beta(e_k) e_k\right). \]
	Let us denote the point
	\[ v=\sum_{k=1}^r  2\beta(e_k)\e_k  e_k\in\partial \W_1, \]
	and define the function $\text{sign}$ by $\text{sign}(x)=x/|x|$ if $x\neq 0$ and $0$ otherwise. For any $\ell'\le \ell_\phi$ we  rewrite $i$ as follows
	\begin{equation*}
		\begin{split}
			i &=\ell'  \left(v+\sum_{k=r+1}^m 2 \beta(e_k) \text{sign}(\lambda_k) e_k \right) + (\ell-\ell')\left( v + \sum_{k=r+1}^m 2\beta(e_k)\left( \frac {\ell}{\ell-\ell'}\lambda_k-\frac {\ell'}{\ell-\ell'}\text{sign}(\lambda_k) \right)  e_k\right)\\
			&=:\ell' w + (\ell-\ell')\left( v +  \sum_{k=r+1}^m 2\beta(e_k)  \lambda'_k e_k\right) .
		\end{split}
	\end{equation*}
	Notice that, since   $\ell\ge 2\ell'$ and $|\lambda_k|\le 1$ it holds $ \left\lvert \lambda'_k \right\rvert \le 1$, thus by formula \eqref{formula facet} we get
	\[v + \sum_{k=r+1}^m 2 \beta(e_k) \lambda'_k  e_k\in \partial \W_1\]
	and therefore $\po(i-\ell' w)=\ell-\ell'$.  We conclude noting that by the hypothesis \eqref{rational zonotope} we can choose $\ell'\le \ell_1$ so that $\ell' w\in\Z^N$, which implies \eqref{additive sd} since $\po(\ell'w)=\ell'.$

	We now prove the last assertion. Since $\po(i)\ge 2\ell_1$, by the previous result there exists $j_0\in (\W_{\ell_1}\setminus \{0\})$ so that $\po(i)=\po(j_0)+\po(i-j_0)$.  Now, if $R-2\ell_1\le \po(j_0)$ we conclude. If not, then $\po(i-j_0)\ge 2\ell_1$ by \eqref{additive sd}, and thus we can find $k_0\in (\W_{\ell_1}\setminus \{0\})$ so that
	\begin{equation}\label{induct0}
	    \po(i-j_0)=\po(k_0)+\po(i-j_0-k_0).
	\end{equation}
	Denoting $j_1=j_0+k_0$, on one hand \eqref{induct0} implies 
	\begin{equation}\label{induct1}
	    \po(i)= \po(j_0)+ \po(j_1-j_0)+\po(i-j_1)\ge \po(j_1)+ \po(i-j_1)
	\end{equation}
    thus equality holds instead. If $\po(j_1)\ge R-2\ell_1$ we conclude, if not \eqref{induct1} yields $\po(i-j_1)\ge 2\ell_1$ and we can iterate. 
    Recalling that  $\po\ge c_\phi>0$ on $\ez\setminus\{0\}$, it is clear  that after a finite number of iterations the process stops, and one can check that  the required properties are satisfied. 
\end{proof}

By the previous lemma it is easy to prove the following result.
\begin{corollary}\label{cor neigh}
	Let  $u\in X$ be a $(1,\phi)$-Lipschitz function and $\ell_1$ as in Lemma \ref{lemma neigh}. Then, for all $ i\in\ez$ it holds
	\[  \sup_{j\in\{ u\ge 0 \}} \left\lbrace u_j-\po(i-j)\right\rbrace  =\max_{j\in\{ u\ge 0 \}}  \left\lbrace u_j-\po(i-j)\right\rbrace.\]
	In addition, if $i\in\{u\le 0\}$, the maximum is reached in a point in $(\{u\le 0\}+\W_{2\e\ell_1})\cap\ez$.
\end{corollary} 

\begin{proof}
 	It is enough to consider $i\in \{u< 0\}\cap \ez$. Let us denote $F= (\{u\le 0\}+\W_{2\e\ell_1})\cap \{ u>0 \}. $  Firstly, by a variant of the argument by iteration employed in the proof of Lemma \ref{lemma neigh}, one can prove that
 	\begin{equation}\label{sup=sup'}
 	    \sup_{j\in\{ u\ge 0 \}} \left\lbrace  u_j-\po(i-j)\right\rbrace =\sup_{j\in F } \left\lbrace u_j-\po(i-j) \right\rbrace.
 	\end{equation}
 	On the other hand, take a point $j_0\in \{u> 0\}$. If $j\in F$ satisfies $u_j-\po(i-j)\ge u_{j_0}-\po(i-{j_0})$, since  $u\le 2\e\ell_1$ in $F$  (as $u$ is $(1,\po)$-Lipschitz) we obtain
 	\[ 2\e\ell_1 +\po(i-j_0)\ge \po(i-j), \]
    which implies that the $\sup$ in \eqref{sup=sup'} is indeed a max.
\end{proof}

\printbibliography 
\end{document}